\documentclass[12pt]{article}
\usepackage{amsmath}
\usepackage{amssymb}
\usepackage{amsthm}
\usepackage{amsfonts}
\usepackage{amscd}
\usepackage[mathscr]{eucal}
\newcommand{\Z} {{\mathbb  Z}}
\newcommand{\Q}{{\mathbb  Q}}
\newcommand{\F}{{\mathbb  F}}
\newcommand{\C}{{\mathbb  C}}

\newcommand{\R} {{\mathbb R}}

\begin{document}
\parindent  25pt
\baselineskip  10mm
\textwidth  15cm    \textheight  23cm \evensidemargin -0.06cm
\oddsidemargin -0.01cm

\title{ {The structures of Hausdorff metric in non-Archimedean spaces}}
\author{\mbox{}
{ Derong Qiu }
\thanks{ \quad E-mail:
derong@mail.cnu.edu.cn, \ derongqiu@gmail.com} \\
(School of Mathematical Sciences, Capital Normal University, \\
Beijing 100048, P.R.China) }

\date{}
\maketitle
\parindent  24pt
\baselineskip  10mm
\parskip  0pt

\par     \vskip  0.6 cm

{\bf Abstract} \ For non-Archimedean spaces $ X $ and $ Y, $ let $
\mathcal{M}_{\flat } (X), \mathfrak{M}(V \rightarrow W) $ and $
\mathfrak{D}_{\flat }(X, Y) $ be the ballean of $ X $ (the family of
the balls in $ X $), the space of mappings from $ X $ to $ Y, $ and
the space of mappings from the ballen of $ X $ to $ Y, $
respectively. By studying explicitly the Hausdorff metric structures
related to these spaces, we construct several families of new metric
structures (e.g., $ \widehat{\rho } _{u}, \widehat{\beta }_{X,
Y}^{\lambda }, \widehat{\beta }_{X, Y}^{\ast \lambda } $) on the
corresponding spaces, and study their convergence, structural
relation, law of variation in the variable $ \lambda , $ including
some normed algebra structure. To some extent, the class $
\widehat{\beta }_{X, Y}^{\lambda } $ is a counterpart of the usual
Levy-Prohorov metric in the probability measure spaces, but it
behaves very differently, and is interesting in itself. Moreover,
when $ X $ is compact and $ Y = K $ is a complete non-Archimedean
field, we construct and study a Dudly type metric of the space of $
K-$valued measures on $ X. $

\par  \vskip  0.4 cm

\par  \vskip  0.4 cm

\par     \vskip  1.2 cm

\hspace{-0.6cm}{\bf 1. \ Introduction }

\par  \vskip  0.2 cm

For non-Archimedean spaces $ X $ and $ Y, $ let $ \mathcal{M}_{\flat
} (X), \mathfrak{M}(V \rightarrow W) $ and $ \mathfrak{D}_{\flat
}(X, Y) $ be the ballean of $ X $ (the family of the balls in $ X
$), the space of mappings from $ X $ to $ Y, $ and the space of
mappings from the ballen of $ X $ to $ Y, $ respectively (see
Sections 2 and 4 for more precise meaning). By studying explicitly
the Hausdorff metric structures related to these spaces, we
construct several families of new metric structures (e.g., $
\widehat{\rho } _{u}, \widehat{\beta }_{X, Y}^{\lambda },
\widehat{\beta }_{X, Y}^{\ast \lambda } $) on the corresponding
spaces (see Sections 3 and 4), and study their convergence,
structural relation, law of variation in the variable $ \lambda , $
including some normed algebra structure (see Theorems 3.2, 3.5, 4.2,
4.5, 4.6, 4.9, 4.11 below). To some extent, the non-Archimedean
class $ \widehat{\beta }_{X, Y}^{\lambda } $ is a counterpart of the
usual Levy-Prohorov metric in the probability measure spaces (see
[D, p.394] and [Ra]), but it behaves very differently, and is
interesting in itself. Moreover, when $ X $ is compact and $ Y = K $
is a complete non-Archimedean field, we construct and study a Dudly
type metric of the space of $ K-$valued measures on $ X $ (see Section5, 
Theorems 5.2, 5.3 below).
\par \vskip 0.2 cm

{\bf Notation and terminology.} \ Let $ X $ be a metric space
endowed with the metric $ d. $ For any subsets $ A $ and $ B $ of $
X , $ the diameter of $ A $ is diam$(A) = \sup \{ d ( x, y ) : x, y
\in A \}, $ the distance between $ A $ and $ B $ is dist$( A, B ) =
\inf \{ d ( x, y ) : x \in A , y \in B \}, $ in particular, for $ x
\in X , $ dist$ ( x, A ) = $ dist$ ( \{ x \} , A ). $ \ For a set $
Y \subset X $ and $ \varepsilon > 0, $ its $ \varepsilon
-$neighborhood is the set $ U_{\varepsilon} ( Y ) = \{ x \in X :
\text{dist} ( x , Y ) < \varepsilon \}. $ The Hausdorff distance
between $ A $ and $ B, $ denoted by $ d_{ H } ( A , B ), $ is
defined by
\begin{align*}
d_{ H } ( A , B ) &= \inf \{ \varepsilon > 0 : \ A \subset
U_{\varepsilon} ( B ) \
\text{ and } \ B \subset U_{\varepsilon} ( A ) \} \\
&= \max \{ \sup_{ a \in A } \text{dist} ( a , B ) , \quad  \sup_{ b
\in B } \text{dist}( b , A ) \} .
\end{align*}
Moreover, for $ a \in X $ and $ r \in (0, \infty ) \subset \R $ (the
real number field), the $ open $ ball of radius $ r $ with center $
a $ is the set $ B_{a} (r) = \{ x \in X : \ d(a, x) < r \}; $ \ The
$ closed $ ball of radius $ r $ with center $ a $ is the set $
\overline{B}_{a} (r) = \{ x \in X : \ d(a, x) \leq r \}. $ \ A ball
in $ (X, d ) $ is a set of the form $ B_{a} (r) $ or $
\overline{B}_{a} (r) $ for some $ a \in X $ and $ r \in (0, \infty
). $ In this paper, a metric may take the value $ \infty . $ For more
basic facts of metric spaces, we refer to [BBI] and [G]. \\
Throughout this paper, unless otherwise specified, we use the
notation $ \widehat{d} $ denote the non-Archimedean metric on a
space $ X , $ that is, $ \widehat{d} $ is a metric on $ X $ and
satisfies the strong triangle inequality \ $  \widehat{d} ( x , z )
\leq \max \{ \widehat{d} ( x , y ), \  \widehat{d} ( y , z ) \} \ (
x , y , z \in X ). $ \ A set $ X $ endowed with a non-Archimedean
metric $ \widehat{d}_{X} $ is called a non-Archimedean metric space
(also called ultrametric space), which is denoted by $ ( X ,
\widehat{d}_{X} ). $ The Hausdorff distance between subsets of $ X $
is then denoted by $ \widehat{d}_{X, H }. $ \ It is well known that,
in $ (X, \widehat{d}_{X} ), $ each ball is both open and closed,
each point of a ball may serve as a center, and a ball may have
infinitely many radii (see [Sc, p.48]). For more basic properties of
non-Archimedean metric spaces, we refer to [BGR], [Sc] and [Q]. \\
As usual, the symbols $ \Z, \Q, \R, \C, \Z_{p}, \Q_{p} $ and $
\F_{p} $ represent the sets of integers, rational numbers, real
numbers, complex numbers, $ p-$adic integers, $ p-$adic numbers and
the field with $ p $ elements, respectively. We denote the
completion of the algebraic closure $ \overline{ \Q_{p}} $ of $
\Q_{p} $ by $ \C_{p}, $ which is endowed the non-Archimedean metric
$ \widehat{d}_{p} $ induced by the normalized valuation $ |.|_{p} $
satisfying $ |p |_{p} = \frac{1}{p} \ (i.e., \widehat{d}_{p} (\alpha
, \beta ) = \mid \alpha - \beta \mid _{p}, \forall \alpha , \beta
\in \C_{p}), $ and called the Tate field (see [K], [Se]).
\par     \vskip 0.3 cm

\hspace{-0.6cm}{\bf 2. \ Ultrametric structures on balls }

\par     \vskip 0.2 cm

Firstly we have the following simple formula for computing the
Hausdorff distance between balls in a non-Archimedean metric space.
\par  \vskip 0.2 cm

{\bf Lemma 2.1.} \ Let $ B_{1} $ and $ B_{2} $ be two balls in a
non-Archimedean metric space $ (X, \widehat{d} ), $ then their
Hausdorff distance is $$ \widehat{d}_{H}( B_{1}, B_{2}) = \left\{
\begin{array}{l} \text{dist}( B_{1}, B_{2})
\quad \text{if} \ B_{1} \bigcap B_{2} = \emptyset ; \\
\max \{ \text{diam} ( B_{1}), \ \text{diam} ( B_{2}) \} \quad
\text{if} \ B_{1} \bigcap B_{2} \neq \emptyset \ \text{and} \ B_{1}
\neq B_{2}; \\
0 \quad \text{if} \ B_{1} = B_{2}.
\end{array}
\right. $$ In particular, $ \widehat{d}_{H}( B_{1}, B_{2}) \geq \max
\{ \text{diam} ( B_{1}), \ \text{diam} ( B_{2}) \} $ if $ B_{1} \neq
B_{2}. $ \\
Moreover, for any non-empty subsets $ A_{i} $ of $ B_{i} \ ( i = 1,
2), $ if $ B_{1} \bigcap B_{2} = \emptyset , $ then $
\widehat{d}_{H} ( A_{1}, A_{2} ) = \widehat{d}_{H}( B_{1}, B_{2}). $
\par  \vskip 0.1 cm

{\bf Proof.} \ Easily follows from the definition.
\par  \vskip 0.2 cm

{\bf Proposition 2.2.} \ Let $ (X, d) $ be a metric space. Then the
metric $ d $ is non-Archimedean if and only if the Hausdorff
distance $ d_{H} $ satisfies the following condition: \ For any
balls $ B_{1} $ and $ B_{2} $ in $ X $ with $ B_{1} \bigcap B_{2} =
\emptyset , $ we have $ d_{H}( B_{1}, B_{2}) = \text{dist}( B_{1},
B_{2}). $
\par  \vskip 0.1 cm

{\bf Proof.} \ The necessity follows from the above Lemma 2.1. For
the sufficiency, we only need to verify that, under the given
condition, the metric $ d $ satisfies the strong triangle
inequality. To see this, let $ a, b, c \in X. $ If $ d (a, b) > d (
a, c), $ then we want to show that $ d (a, b) = d (c, b). $ In fact,
one can take a real number $ r $ with $ d ( a, c) < r < d (a, b). $
Then for the $ closed $ ball $ \overline{B} = \overline{B}_{a} (r),
$ we have $ c \in \overline{B} $ and $ b \notin \overline{B}. $ Take
an $ open $ ball $ B_{1}= B_{b}(r_{1}) $ of radius $ r_{1} > 0 $
with center $ b $ in $ X $ such that $ B_{1} \bigcap \overline{B} =
\emptyset . $ By the given condition, we have $ d_{H}( B_{1},
\overline{B} ) = \text{dist}( B_{1}, \overline{B} ), $ then it
follows that $ \text{dist} ( \overline{b}_{1}, B_{1} ) =
\text{dist}( \overline{b}_{2}, B_{1} ) $ for any $ \overline{b}_{1},
\overline{b}_{2} \in \overline{B}. $ In particular, we have $
\text{dist} ( a, B_{1} ) = \text{dist}( c, B_{1} ). $ Then by taking
the limit $ r_{1} \longrightarrow 0, $ we get $ d (a, b) = d (c, b).
$ So we are done, and the other cases are obvious. The proof of
Proposition 2.2 is completed.  \quad $ \Box $
\par  \vskip 0.2 cm

{\bf Lemma 2.3.} \ Let $ (X, \widehat{d} ) $ be a non-Archimedean
metric space and $ A \subset X $ be a non-empty subset. \\
(1) \ For any $ \varepsilon > 0, $ we have \ $ U_{\varepsilon } (
U_{\varepsilon } ( A ) ) = U_{\varepsilon } ( A ). $ In particular,
if $ A = B $ is a ball with diameter $ r = \text{diam} (B), $ then $
U_{\varepsilon } ( B ) = \left\{ \begin{array}{l} B \quad \text{if}
\ \varepsilon \leq r, \\
B_{b}(\varepsilon ) \quad \text{if} \ \varepsilon > r.
\end{array} \right. $ Here $ b \in B $ can be arbitrarily chosen.
Furthermore, for any positive numbers $ \varepsilon _{1} $ and $
\varepsilon _{2}, $ we have \ $ U_{\varepsilon _{2}} (
U_{\varepsilon _{1}}(B)) = U_{ \max \{\varepsilon _{1}, \
\varepsilon _{2} \}}(B). $ \\
(2) \ Let $ \varepsilon > 0, $ if $ U_{\varepsilon } ( A ) \neq X, $
then the characteristic function $$ \chi _{U_{\varepsilon } ( A
)}(x) = \left\{ \begin{array}{l} 1 \quad \text{if} \ x \in
U_{\varepsilon } ( A ), \\
0 \quad \text{if} \ x \notin U_{\varepsilon } ( A )
\end{array} \right. $$ is uniformly continuous, and
$ \text{dist}(U_{\varepsilon } ( A ), \ X \setminus U_{\varepsilon }
( A ) ) \geq \varepsilon . $
\par  \vskip 0.1 cm

{\bf Proof.} \ (1) follows easily from the definitions. \\
(2) \ By (1), we know that $ U_{\varepsilon } ( A ) $ is uniformly
open, so $ \chi _{U_{\varepsilon } ( A )} $ is uniformly continuous
(see [Sc, p.50]). The last inequality follows easily from the
definition. \quad $ \Box $
\par  \vskip 0.2 cm

{\bf Lemma 2.4.} \ Let $ (X, \widehat{d} ) $ be a non-Archimedean
metric space. Then  \\
(1) \ $ \widehat{d}_{H} $ is a non-Archimedean semi-metric on $
2^{X} \setminus \{ \emptyset \} $ ( the set of all subsets of $ X $
except the empty set $ \emptyset $ ). \\
(2) \ $  \widehat{d}_{H} ( A, \ \overline{A} ) = 0 $ for any $
\emptyset \neq A \subset X, $ where $ \overline{A} $ denote the
closure of $ A $ in $ X. $ \\
(3) \ If $ A $ and $ B $ are closed subsets of $ X $ and $
\widehat{d}_{H} ( A, \ B ) = 0, $ then $ A = B. $
\par  \vskip 0.1 cm

{\bf Proof.} \ By Proposition7.3.3 of [BBI, p.252], we only need
verify the strong triangle inequality for $ \widehat{d}_{H}. $ To
see this, take $ A, B, C \in 2^{X} \setminus \{ \emptyset \}, $ we
want to show that $ \widehat{d}_{H} ( A, \ B ) \leq \max \{
\widehat{d}_{H} ( A, \ C ), \ \widehat{d}_{H} ( C, \ B ) \}. $ We
may as well assume that $ \widehat{d}_{H} ( A, \ C ) \leq
\widehat{d}_{H} ( C, \ B ) = r. $ If $ r = \infty , $ then we are
done. So we assume that $ r < \infty . $ Then for any $ \varepsilon
> 0, $  by definition and Lemma 2.3(1), it is easy to see that \\
$ A \subset U_{r + \varepsilon } (C) \subset U_{r + \varepsilon }
(U_{r + \varepsilon } (B)) = U_{r + \varepsilon } (B), \\
B \subset U_{r + \varepsilon } (C) \subset U_{r + \varepsilon }
(U_{r + \varepsilon } (A)) = U_{r + \varepsilon } (A). $ \\
So $ \widehat{d}_{H} ( A, \ B ) \leq r + \varepsilon $ for any $
\varepsilon > 0. $ This shows that \\
$ \widehat{d}_{H} ( A, \ B )
\leq r = \max \{ \widehat{d}_{H} ( A, \ C ), \ \widehat{d}_{H} ( C,
\ B ) \}. $ \quad $ \Box $ \\
It follows easily by Lemma 2.4 that $ (\mathcal{M}(X), \
\widehat{d}_{H} ) $ is a non-Archimedean metric space, where $
\mathcal{M}(X) $ is the set of non-empty closed subsets of $ X. $
Moreover, every element of the quotient $ ( 2^{X} \setminus
\{\emptyset \}) / \widehat{d}_{H} $ can be represented by a closed
set and therefore $ ( 2^{X} \setminus \{\emptyset \}) /
\widehat{d}_{H} $ is naturally identified with $ ( \mathcal{M}(X), \
\widehat{d}_{H} ). $ Furthermore, if $ X $ is complete (resp.
compact), then $ ( \mathcal{M}(X), \ \widehat{d}_{H} ) $ is complete
(resp. compact) (see [BBI, Prop.7.3.7 and Thm. 7.3.8]). \\
Recall that a ball in $ (X, \widehat{d}) $ is a set of the form
$ B_{a} (r) $ or $ \overline{B}_{a} (r) $ for some $ a \in X $ and
$ r \in (0, \infty ), $ and every ball in a non-Archimedean metric
space is both open and closed.
\par  \vskip 0.2 cm

{\bf Definition 2.5.} \ Given a non-Archimedean metric space $ (X, \
\widehat{d}), $ we set \\
$ \mathcal{M}_{\flat } (X) = \{ \text{ all balls in } X \}, \quad
\overline{\mathcal{M}}_{\flat } (X) = \mathcal{M}_{\flat } (X)
\bigcup \{ \{x \} : \ x \in X \}. $ \\
Since each ball in $ X $ is both open and closed, it is obvious that
$ \mathcal{M}_{\flat } (X), \ \overline{\mathcal{M}}_{\flat } (X)
\subset \mathcal{M}(X). $ So both $ ( \mathcal{M}_{\flat } (X), \
\widehat{d}_{H} ) $ and $ ( \overline{\mathcal{M}}_{\flat } (X), \
\widehat{d}_{H} ) $ are subspaces of $ (\mathcal{M}(X), \
\widehat{d}_{H} ), $ and they are all non-Archimedean metric spaces.
\par  \vskip 0.2 cm

{\bf Lemma 2.6.} \ Let $ (X, \widehat{d} ) $ be a non-Archimedean
metric space. If $ \{ B_{m} \}_{m = 1}^{\infty } \subset (
\mathcal{M}_{\flat } (X), \ \widehat{d}_{H} ) $ is a Cauchy sequence
of balls, then either $ \text{diam} ( B_{m} ) \rightarrow 0 $ as $ m
\rightarrow \infty , $ or $ B_{m} = B_{m + i } \ ( i = 1, 2, \cdots
) $ for sufficiently large $ m. $
\par  \vskip 0.1 cm

{\bf Proof.} \ If the conclusion that $ \text{diam} ( B_{m} )
\rightarrow 0 $ as $ m \rightarrow \infty $ does not hold, then
there exists an $ \varepsilon _{0} > 0 $ such that for any $ N
> 0, $ there exists an integer $ n > N $ satisfying $ \text{diam} (
B_{n} ) \geq \varepsilon _{0}. $ Put $ \varepsilon =
\frac{\varepsilon _{0}}{2}. $ Since $ \{ B_{m} \}_{m = 1}^{\infty }
$ is a Cauchy sequence, there exists a positive integer $ m_{0} $
such that for any positive integers $ m, n \geq m_{0}, $ we have $
\widehat{d}_{H} ( B_{m}, \ B_{n} ) < \varepsilon . $ Take an integer
$ n_{0} > m_{0} $ such that $ \text{diam} ( B_{n_{0}} ) \geq
\varepsilon _{0}, $ then for any positive integer $ i, $ we must
have $ B_{n_{0}} = B_{n_{0} + i}. $
Otherwise, by Lemma 2.1, one has \\
$ \varepsilon _{0} \leq \max \{\text{diam} ( B_{n_{0}} ), \
\text{diam} ( B_{n_{0} + i} ) \} \leq \widehat{d}_{H}(B_{n_{0}}, \
B_{n_{0} + i} ) < \varepsilon  = \frac{\varepsilon _{0}}{2}, $ \\
a contradiction! So we obtain that $ B_{n_{0}} = B_{n_{0} + i } $
for all positive integers $ i. $ This completes the proof of
Lemma 2.6. \quad $ \Box $ \\
If the sequence $ \{ A_{i} \}_{i = 1}^{\infty } $ of subsets of $
(X, \widehat{d}) $ is convergent to a subset $ A $ under $
\widehat{d}_{H}, $ we will write it as \ $ \widehat{d}_{H}- \lim _{
i \rightarrow \infty } A_{i} = A. $
\par  \vskip 0.2 cm

{\bf Proposition 2.7.} \ Let $ (X, \widehat{d} ) $ be a
non-Archimedean metric space and $ \{ B_{m} \}_{m = 1}^{\infty }
\subset ( \mathcal{M}_{\flat } (X), \ \widehat{d}_{H} ) $ be a
sequence of balls. If \ $ \widehat{d}_{H}- \lim _{ m \rightarrow
\infty } B_{m} = D  $ for a non-empty closed subset $ D $ of $ X, $
then $ D \in \overline{\mathcal{M}}_{\flat } (X). $
\par  \vskip 0.1 cm

{\bf Proof.} \ Since $ \{ B_{m} \}_{m = 1}^{\infty } $ is a Cauchy
sequence, by Lemma 2.6 above, $ B_{m} = B_{m + i} \ ( i = 1, 2,
\cdots ) $ for sufficiently large $ m, $ or $ \text{diam} ( B_{m} )
\rightarrow 0 $ as $ m \rightarrow \infty . $ For the former, we
have $ D = B_{m} \in \mathcal{M}_{\flat } (X). $ So we only need to
verify that $ D $ is an one-point set if
$ \text{diam} (B_{m}) \rightarrow 0 $ as $ m \rightarrow \infty . $
In fact, if otherwise, then $ \sharp D > 1 $ and so there exist
elements $ a, b \in D $ with $ a \neq b. $ Denote $ r = \widehat{d}
(a, b), $ then $ r > 0. $ Take $ \varepsilon = \frac{r}{2}, $ then
there exists a positive integer $ m_{0} $ such that for any positive
integer $ m > m_{0} $ we have both $ \text{diam} (B_{m}) < \varepsilon $
and $ \widehat{d}_{H} ( B_{m}, \ D ) < \varepsilon . $ Thus by
definition we get $ \text{dist} ( a, \ B_{m}) \leq \widehat{d}_{H} (
B_{m}, \ D ) < \varepsilon $ \ and $ \text{dist} ( b, \ B_{m}) \leq
\widehat{d}_{H} ( B_{m}, \ D ) < \varepsilon . $ Then for any $ c
\in B_{m}, $ since
$ \left\{
\begin{array}{l} \widehat{d} ( a, c ) = \text{dist} ( a, \ B_{m}) \
\text{if} \ a \notin B_{m};  \\
\widehat{d} ( a, c ) \leq \text{diam} (B_{m}) \ \text{if} \ a \in B_{m},
\end{array}
\right. $
we get $ \widehat{d} ( a, c ) < \varepsilon . $ Likewise, $
\widehat{d} ( b, c ) < \varepsilon . $ Hence by the strong triangle
inequality, we get $ \widehat{d} ( a, b ) < \varepsilon . $ This is
impossible because $ \widehat{d} ( a, b ) = r = 2 \varepsilon . $
Therefore $ \sharp D = 1 $ and so we are done. This completes the
proof of Proposition 2.7. \quad $ \Box $
\par  \vskip 0.2 cm

{\bf Theorem 2.8.} \ Let $ (X, \widehat{d} ) $ be a non-Archimedean
metric space. \\
(1) \ If $ (X, \widehat{d}) $ is complete, then $ (
\overline{\mathcal{M}}_{\flat } (X), \ \widehat{d}_{H} ) $ is
complete, and $ \mathcal{M}_{\flat } (X) $ is dense in $
\overline{\mathcal{M}}_{\flat } (X), $ i.e., $
\overline{\mathcal{M}}_{\flat } (X) $ is the completion of $
\mathcal{M}_{\flat } (X). $ \\
(2) \ If $ (X, \widehat{d}) $ is compact, then $ (
\overline{\mathcal{M}}_{\flat } (X), \ \widehat{d}_{H} ) $ is
compact.
\par  \vskip 0.1 cm

{\bf Proof.} \ (1) \ To show $ \overline{\mathcal{M}}_{\flat } (X) $
is complete, we need to verify that every Cauchy sequence in it has
a limit. So we take a Cauchy sequence $ \{C_{n} \}_{n = 1}^{\infty }
$ in $ ( \overline{\mathcal{M}}_{\flat } (X), \ \widehat{d}_{H} ). $
Then it has a limit $ C $ in $ \mathcal{M} ( X ) $ because $ (
\mathcal{M} ( X ), \ \widehat{d}_{H} ) $ is complete by our
assumption. We need to show that $ C \in
\overline{\mathcal{M}}_{\flat } (X). $ To see this, firstly we
assume that there exists a positive integer $ n_{0} $ such that $
C_{n} \in \mathcal{M}_{\flat } (X) $ for all positive integers $ n >
n_{0}, $ then by Proposition 2.7 we get $ C \in
\overline{\mathcal{M}}_{\flat } (X). $ Next we assume that there
exists a subsequence $ \{C_{n_{i}} \}_{i = 1}^{\infty } $ of $
\{C_{n} \}_{n = 1}^{\infty } $ such that $ C_{n_{i}} = \{ x_{i} \} $
with $ x_{i} \in X $ for all $ i = 1, 2, \cdots . $ Obviously, $
\{C_{n_{i}} \}_{i = 1}^{\infty } $ is convergent and $
\widehat{d}_{H}- \lim _{ i \rightarrow \infty } C_{n_{i}} = C. $ As
one can easily see that $ \widehat{d}_{H} ( \{ x_{i} \}, \ \{ x_{j}
\}) = \widehat{d}( x_{i} , \ x_{j} ), $ it follows that $ \{x_{i}
\}_{i = 1}^{\infty } $ is a Cauchy sequence in $ ( X, \widehat{d}),
$ and so $ \lim _{ i \rightarrow \infty } x_{i} = x $ for some $ x
\in X. $ Hence $ \widehat{d}_{H} ( \{ x_{i} \}, \ \{ x \}) =
\widehat{d}( x_{i} , \ x ) \rightarrow 0  $ as $ i \rightarrow
\infty . $ So $ \widehat{d}_{H}- \lim _{ i \rightarrow \infty }
C_{n_{i}} = \{ x \}, $ therefore $ C = \{ x \} \in
\overline{\mathcal{M}}_{\flat } (X). $ This shows that $ (
\overline{\mathcal{M}}_{\flat } (X), \ \widehat{d}_{H} ) $ is
complete. In particular, $ \overline{\mathcal{M}}_{\flat } (X) $ is
a closed subset in $ \mathcal{M} ( X ). $ \\
Now we come to show the density of $ \mathcal{M}_{\flat } (X) $ in $
\overline{\mathcal{M}}_{\flat } (X). $ We may as well assume that $
\mathcal{M}_{\flat } (X) \neq \overline{\mathcal{M}}_{\flat } (X). $
Let $ C \in \overline{\mathcal{M}}_{\flat } (X) \setminus
\mathcal{M}_{\flat } (X). $ Then $ C = \{ x \} $ for some $ x \in X.
$ For every positive integer $ n, $ let $ C_{n} = B_{x} (
\frac{1}{n}) $ be the $ open $ ball of radius $ \frac{1}{n} $ with
center $ x $ as defined before. Then $ C_{n} \in \mathcal{M}_{\flat
} (X) $ for all $ n \geq 1. $ By Lemma 2.1, we have $
\widehat{d}_{H} ( C, \ C_{n} ) = \widehat{d}_{H} ( \{ x \}, \ C_{n}
) = \text{diam} ( B_{x} ( \frac{1}{n}) ) \leq \frac{1}{n}
\rightarrow 0 $ as $ n \rightarrow \infty . $ Therefore $
\widehat{d}_{H}- \lim _{ n
\rightarrow \infty } C_{n} = C. $ This proves the density. \\
(2) \ As mentioned above, $ ( \mathcal{M} (X), \ \widehat{d}_{H} ) $
is compact by our assumption, then the conclusion follows easily
from (1). This completes the proof of Theorem 2.8. \quad $ \Box $
\par  \vskip 0.2 cm

{\bf Remark 2.9.} \ Let $ (X, \widehat{d} ) $ be a complete
non-Archimedean metric space. By Theorem 2.8, it is easy to see that
$ ( \mathcal{M}_{\flat } (X), \ \widehat{d}_{H} ) $ is complete if
and only if all the one-point subsets $ \{ x \} $ of $ X $ are
balls, in other words, $ (X, \widehat{d} ) $ is discrete.
\par  \vskip 0.2 cm

{\bf Proposition 2.10.} \ The function $ \rho _{f}: \
\overline{\mathcal{M}}_{\flat } (X) \rightarrow \R, $
$$  \rho _{f} ( A ) = f ( \text{diam} (A)), \quad  A \in
\overline{\mathcal{M}}_{\flat } (X), $$ is
uniformly continuous for every $ f : \ \R_{ \geq 0 } \rightarrow \R $
which is continuous at the point $ 0. $
\par  \vskip 0.1 cm

{\bf Proof.} \ For any $ \varepsilon > 0, $ by the continuity, $
\exists \delta > 0 $ such that $ \mid f ( x ) - f ( 0 ) \mid
_{\infty } < \frac{ \varepsilon }{2} $ for any real number $ x \in [
0, \ \delta ), $ where $ \mid \ \mid _{\infty } $ denotes the usual
absolute value of the real number field $ \R. $ Then for any $
A_{1}, \ A_{2} \in \overline{\mathcal{M}}_{\flat } (X) $ satisfying
$ \widehat{d}_{H} ( A_{1}, \ A_{2} ) < \delta , $ by Lemma 2.1 one
has either $ A_{1} = A_{2} $ or $ \max \{ \text{diam} (A_{1}), \
\text{diam} (A_{2}) \} \leq \widehat{d}_{H} ( A_{1}, \ A_{2} ) <
\delta . $ So either \ $ \rho _{f} ( A_{1} ) = \rho _{f} ( A_{2} ) $
\ or \ $ \mid \rho _{f} ( A_{1} ) - \rho _{f} ( A_{2} ) \mid
_{\infty } = \mid f ( \text{diam} (A_{1}))  - f ( \text{diam}
(A_{2}))  \mid _{\infty } \leq \ \mid f ( \text{diam} (A_{1})) - f (
0 ) \mid _{\infty } + \mid f ( \text{diam} (A_{2})) - f ( 0 ) \mid
_{\infty } < \frac{\varepsilon}{2} + \frac{\varepsilon}{2} =
\varepsilon . $ This shows that $ \rho _{f} $ is uniformly
continuous, and the proof of Proposition 2.10 is completed. \quad $
\Box $
\par  \vskip 0.2 cm

{\bf Definition 2.11.} \ Given a non-Archimedean metric space $ (X,
\widehat{d} ), $ in the above Definition 2.5 we have defined two
corresponding non-Archimedean metric spaces $ ( \mathcal{M}_{\flat
}(X), \ \widehat{d}_{H}) $ and $ ( \overline{\mathcal{M}}_{\flat
}(X), \ \widehat{d}_{H} ). $  Now we view the symbols $
\mathcal{M}_{\flat } $ and $ \overline{\mathcal{M}}_{\flat } $ as
operators, and define inductively two sequences of metric spaces $
\mathcal{M}_{\flat }^{(n)} (X) $ and
$ \overline{\mathcal{M}}_{\flat }^{(n)} (X) $ as follows: \\
Set $ \mathcal{M}_{\flat }^{(1)} (X) =  \mathcal{M}_{\flat }(X);
\quad \overline{\mathcal{M}}_{\flat }^{(1)} (X) =
\overline{\mathcal{M}}_{\flat } (X) $ \\
with the induced Hausdorff distance $  \widehat{d}_{H}^{(1)}
= \widehat{d}_{H}; $ \\
$ \mathcal{M}_{\flat }^{(2)} (X) = \mathcal{M}_{\flat } (
\mathcal{M}_{\flat }^{(1)}(X)); \quad \overline{\mathcal{M}}_{\flat
}^{(2)} (X) = \overline{\mathcal{M}}_{\flat } (
\overline{\mathcal{M}}_{\flat }^{(1)}(X) ) $ \\
with the induced Hausdorff distance $
\widehat{d}_{H}^{(2)} = (\widehat{d}_{H}^{(1)})_{H}; $ \\
And for general positive integer $ n > 1, $ set \\
$ \mathcal{M}_{\flat }^{(n)} (X) = \mathcal{M}_{\flat } (
\mathcal{M}_{\flat }^{(n - 1)}(X)); \quad
\overline{\mathcal{M}}_{\flat }^{(n)} (X) =
\overline{\mathcal{M}}_{\flat } (
\overline{\mathcal{M}}_{\flat }^{(n - 1)}(X) ) $ \\
with the induced Hausdorff distance $ \widehat{d}_{H}^{(n)} = (
\widehat{d}_{H}^{(n - 1)})_{H}. $ \\
For convenience, we write $ \overline{\mathcal{M}}_{\flat }^{(0)}(X)
= X = \mathcal{M}_{\flat }^{(0)} (X). $ It is easy to know that all
$ ( \mathcal{M}_{\flat }^{(n)} (X), \ \widehat{d}_{H}^{(n)} ) $ and
$ ( \overline{\mathcal{M}}_{\flat }^{(n)} (X), \
\widehat{d}_{H}^{(n)} ) $ are non-Archimedean metric spaces.
In particular, if $ ( X, \widehat{d}) $ is complete (resp. compact),
then all $ ( \overline{\mathcal{M}}_{\flat }^{(n)} (X), \
\widehat{d}_{H}^{(n)} ) $ are complete (resp. compact). \\
If we set $ j_{X} : \ X \rightarrow \overline{\mathcal{M}}_{\flat
}(X), \ x \mapsto \{ x \}, $ \ then obviously $ \widehat{d}_{H}
(j_{X}(x_{1}), j_{X}(x_{2})) = \widehat{d} (x_{1}, x_{2}), $ so $
j_{X} $ is an isometric embedding and $ X $ can be viewed as a
subspace of $ \overline{\mathcal{M}}_{\flat }(X). $ By this way, for
each non-negative integer $ n,  \overline{\mathcal{M}}_{\flat
}^{(n)} (X) $ can be viewed as a subspace of $
\overline{\mathcal{M}}_{\flat }^{(n + 1)} (X). $ So we may define
the direct limit space $ \overline{\mathcal{M}}_{\flat }^{(\infty )}
(X) = \lim _{n \rightarrow \infty } \overline{\mathcal{M}}_{\flat
}^{(n)} (X) $ endowed with the inductive topology. On the other hand,
it is easy to see that $ \overline{\mathcal{M}}_{\flat }^{(\infty )} (X) $
is a non-Archimedean metric space endowed with the metric $
\widehat{d}_{H}^{(\infty)} $ defined as follows:
For any elements $ \alpha , \beta \in \overline{\mathcal{M}}_{\flat
}^{(\infty )} (X), $ there exists a non-negative integer $ n $ such
that $ \alpha , \beta \in \overline{\mathcal{M}}_{\flat }^{(n)} (X),
$ and we define $ \widehat{d}_{H}^{(\infty )}(\alpha , \beta ) =
\widehat{d}_{H}^{(n)}(\alpha , \beta ). $ Let $
\widehat{\overline{\mathcal{M}}}_{\flat }^{(\infty )} (X) $ be the
completion of $ \overline{\mathcal{M}}_{\flat }^{(\infty )} (X) $
under the metric $ \widehat{d}_{H}^{(\infty )}. $ A question is
\par  \vskip 0.2 cm

{\bf Question 2.12.} \ Is it true that $
\widehat{\overline{\mathcal{M}}}_{\flat }^{(\infty )} (X) =
\overline{\mathcal{M}}_{\flat }^{(\infty )} (X), $ in other words,
is $ \overline{\mathcal{M}}_{\flat }^{(\infty )} (X) $ complete ?

\par     \vskip  0.3 cm

\hspace{-0.6cm}{\bf 3. \ Ultrametric structures on mappings}

\par  \vskip  0.2 cm

For any two non-Archimedean metric spaces $ (X, \widehat{d}_{X} ) $ and
$ (Y, \widehat{d}_{Y} ), $ let $ Z = X \times Y $ be their Cartesian
product. As usual, we define a function $ \widehat{d}_{Z} $ on
$ Z \times Z $ by
$$ \widehat{d}_{Z} ((x_{1}, y_{1}), (x_{2}, y_{2}))
= \max \{\widehat{d}_{X} (x_{1}, x_{2}), \ \widehat{d}_{Y} (y_{1},
y_{2}) \} \ (\forall \ (x_{1}, y_{1}), (x_{2}, y_{2}) \in Z). $$
Then it is easy to see that $ \widehat{d}_{Z} $ is a non-Archimedean
metric on $ Z, $ which is called the square metric of $
\widehat{d}_{X} $ and $ \widehat{d}_{Y} $(see [M]). Throughout this
paper, for any metric spaces $ V $ and $ W, $ we denote
\begin{align*}
&\mathfrak{M}(V \rightarrow W) =
\{\text{all maps} \ f : V \rightarrow W \}, \\
&\mathfrak{C}(V \rightarrow W)
= \{\text{all continuous maps} \ f : V \rightarrow W \}, \\
&\mathfrak{C}_{u}(V \rightarrow W) = \{\text{all uniformly continuous maps} \
f : V \rightarrow W \}.
\end{align*}
Moreover, for every $ f \in \mathfrak{M}(V \rightarrow W), $ we write
$ \Gamma _{f} = \{(v, f(v)) : \ v \in V \} \subset V \times W, $ which
is the graph of the map $ f. $
\par     \vskip  0.2 cm

{\bf Definition 3.1.} \ Let $ (X, \widehat{d}_{X} ) $ and $ (Y, \widehat{d}_{Y} ) $
be two non-Archimedean metric spaces. For any $ f, g \in
\mathfrak{M}(X \rightarrow Y), $ we define
$$ \widehat{\rho} _{H}(f, g) = \widehat{d}_{Z, H} (\Gamma _{f}, \Gamma _{g}),
\quad  \widehat{\rho } _{s}(f, g) = \sup _{x \in X} \widehat{d}_{Y}
(f(x), g(x)), $$ where $ Z = X \times Y $ is endowed the square
metric $ \widehat{d}_{Z} $ of $ \widehat{d}_{X} $ and $
\widehat{d}_{Y} $ defined as above, and $ \widehat{d}_{Z, H} $ is
the corresponding Hausdorff distance between subsets of $ Z. $ We
also define for any real number $ \varepsilon > 0 $
\begin{align*}
&\theta _{f} (\varepsilon ) = \sup \{\widehat{d}_{Y}(f(x_{1}), f(x_{2})) :
\ x_{1}, x_{2} \in X \ \text{and} \ \widehat{d}_{X}(x_{1}, x_{2}) < \varepsilon
\} \ \text{and} \\
&\theta _{f, g} (\varepsilon ) = \min \{ \theta _{f}(\varepsilon ), \
\theta _{g}(\varepsilon ) \}.
\end{align*}
Obviously, for any $ a \geq 0, $ the right-hand limit
$ \lim _{\varepsilon \rightarrow a^{+}} \theta _{f, g} (\varepsilon ) $
exists, and we denote it by $ \theta _{f, g} (a)^{+}. $ From Lemma 2.4 above,
$ \widehat{\rho } _{H} $ is a non-Archimedean semi-metric in
$ \mathfrak{M}(X \rightarrow Y), $ and
$ (\mathfrak{M}(X \rightarrow Y) / \widehat{\rho } _{H}, \widehat{\rho } _{H} ) $
is a non-Archimedean metric space (see [BBI, p.2] for the notation).
\par     \vskip  0.2 cm

Recall that for a map $ f \in \mathfrak{M}(X \rightarrow Y), $ its $ distorsion \
 \text{dis} f $ is defined by $$ \text{dis} f = \sup _{x_{1}, x_{2} \in X}
\mid \widehat{d}_{Y} (f(x_{1}), f(x_{2})) - \widehat{d}_{X} (x_{1}, x_{2}) \mid
\ (\text{see [BBI], p.249}). $$
\par     \vskip  0.2 cm

{\bf Theorem 3.2.} \ Let $ (X, \widehat{d}_{X} ) $ and $ (Y, \widehat{d}_{Y} ) $
be two non-Archimedean metric spaces.  \\
(1) \ For any $ f, g \in \mathfrak{M}(X \rightarrow Y), $ we have
\begin{align*}
&\widehat{\rho } _{H}(f, g) \leq \widehat{\rho } _{s}(f, g)
\leq \max \{ \theta _{f, g} (\widehat{\rho } _{H}(f, g))^{+}, \
\widehat{\rho } _{H}(f, g) \}, \ \text{and} \\
&\theta _{f, g} (\widehat{\rho } _{H}(f, g))^{+} \leq \widehat{\rho } _{H}(f, g)
+ \min \{\text{dis} f, \ \text{dis} g \}.
\end{align*}
(2) \ Let $ f_{n} \in \mathfrak{M}(X \rightarrow Y)
\ (n = 1, 2, \cdots , \infty ) $ and $ g \in \mathfrak{C}_{u}(X \rightarrow Y). $
Then the sequence $ \{ f_{n} \}_{n = 1}^{\infty } $ converge uniformly
to $ g $ if and only if
$ \lim _{n \rightarrow \infty } \widehat{\rho } _{H}(f_{n}, g) = 0. $
\par   \vskip  0.1 cm
{\bf Proof.} \ (1) \ The first inequality follows easily from the definitions.
For the second inequality, let $ \varepsilon > \widehat{\rho } _{H}(f, g), $ then, for
any $ x \in X, $ by definition, $ \varepsilon > \text{dist} ((x, f(x)), \Gamma _{g})
= \inf _{x^{\prime } \in X} \max \{\widehat{d}_{X} (x, x^{\prime }),
\widehat{d}_{Y} (f(x), g(x^{\prime })) \}, $ where $ (x, f(x)) \in \Gamma _{f}. $
So there exists a point $ x^{\prime } \in X $ such that
$ \widehat{d}_{X} (x^{\prime }, x) < \varepsilon $ and
$ \widehat{d}_{Y} (f(x), g(x^{\prime })) < \varepsilon . $ Then
$$ \widehat{d}_{Y} (f(x), g(x)) \leq \max \{\widehat{d}_{Y} (f(x), g(x^{\prime })),
\ \widehat{d}_{Y} (g(x^{\prime }), g(x)) \} \leq \max \{ \varepsilon , \
\theta _{g} (\varepsilon ) \}, $$ which implies $ \widehat{\rho } _{s}(f, g)
\leq \max \{ \varepsilon , \ \theta _{g} (\varepsilon ) \} $ because $ x $ is
arbitrary. Similarly we have $ \widehat{\rho } _{s}(f, g)
\leq \max \{ \varepsilon , \ \theta _{f} (\varepsilon ) \}. $ Therefore we obtain
$$ \widehat{\rho } _{s}(f, g) \leq \max \{ \varepsilon ,
\min \{ \theta _{f} (\varepsilon ),
\theta _{g} (\varepsilon ) \} \} = \max \{ \varepsilon ,
\theta _{f, g} (\varepsilon ) \} $$ for arbitrary $ \varepsilon >
\widehat{\rho } _{H}(f, g), $ and the second inequality follows.
For the last inequality, let $ \varepsilon > \widehat{\rho } _{H}(f, g), $
for any $ x_{1}, x_{2} \in X $ with $ \widehat{d}_{X} (x_{1}, x_{2})
< \varepsilon , $ by the definition of $ \text{dis} f, $ we have
$ \widehat{d}_{Y} (f(x_{1}), f(x_{2})) <
\varepsilon + \text{dis} f. $ So $ \theta _{f} (\varepsilon ) \leq
\varepsilon + \text{dis} f. $ Similarly, $ \theta _{g} (\varepsilon )
\leq \varepsilon + \text{dis} g. $ Therefore $ \theta _{f, g} (\varepsilon )
\leq \varepsilon + \min \{\text{dis} f, \ \text{dis} g \} $ for all
$ \varepsilon > \widehat{\rho } _{H}(f, g), $ by taking the right-hand
limit, the last inequality follows. This proves (1). \\
(2) \ If $ \{ f_{n} \}_{n = 1}^{\infty } $ converge uniformly to $ g, $
then by definition,
$$ \widehat{\rho } _{s}(f_{n}, g) = \sup _{x \in X} \widehat{d}_{Y} (f_{n}(x), g(x))
\rightarrow 0 \ \text{as} \ n \rightarrow \infty . $$ So by (1) above,
$ \widehat{\rho } _{H}(f_{n}, g) \rightarrow 0 $ as $ n \rightarrow \infty . $ \\
Conversely, assume that $ \lim _{n \rightarrow \infty } \widehat{\rho } _{H}(f_{n}, g)
= 0. $ Since $ g $ is uniformly continuous, for any $ \varepsilon > 0, $ there
exists a $ \delta > 0 $ such that $ \widehat{d}_{Y} (g(x_{1}), g(x_{2})) <
\varepsilon $ for all $ x_{1}, x_{2} \in X $ with $ \widehat{d}_{X}
(x_{1}, x_{2}) < \delta . $ So by definition, $ \theta _{g}(\delta ) \leq
\varepsilon . $ Denote $ \varepsilon _{0} = \min \{ \varepsilon , \delta \}, $
then by assumption, there exists a positive integer $ N $ such that
$ \widehat{\rho } _{H}(f_{n}, g) < \varepsilon _{0} $ for all $ n > N. $ Obviously,
$ \theta _{f_{n}, g} (\eta ) \leq \theta _{g} (\eta ) \leq \theta _{g} (\delta )
\leq \varepsilon $ for all $ \eta : \widehat{\rho } _{H}(f_{n}, g) < \eta \leq
\varepsilon _{0}, $ which implies that
$ \theta _{f_{n}, g}(\widehat{\rho } _{H}(f_{n}, g))^{+} \leq \varepsilon . $
Therefore, by (1) above, $ \widehat{\rho } _{s}(f_{n}, g) \leq \varepsilon $
for all $ n > N, $ which implies that $ \{ f_{n} \}_{n = 1}^{\infty } $ converge
uniformly to $ g, $ this proves (2), and the proof of Theorem 3.2 is completed.
\quad $ \Box $
\par     \vskip  0.2 cm

{\bf Remark 3.3.} \ For $ f, g $ as in Theorem 3.2(1) above, if one of them
is a nonexpanding map (see [BBI], p.9), then it can be easily verified that
$ \widehat{\rho } _{H}(f, g) = \widehat{\rho } _{s}(f, g). $
\par   \vskip  0.2 cm

{\bf Definition 3.4.} \ Let $ (X, \widehat{d}_{X} ) $ and $ (Y, \widehat{d}_{Y} ) $
be two non-Archimedean metric spaces. \\
(1) \ For any $ f, g \in
\mathfrak{M}(X \rightarrow Y), $ we define
\begin{align*}
\widehat{\rho } _{b}(f, g) = \inf \{\varepsilon > 0 : \
\widehat{d}_{Y} (f(x), g(x^{\prime }))
< \varepsilon \ \text{and} \ \widehat{d}_{Y} (g(x), f(x^{\prime \prime }))
< \varepsilon \\
\text{for some} \ x^{\prime }, x^{\prime \prime }
\in B_{x}(\varepsilon ) \ (\forall x \in X) \}.
\end{align*}
It is easy to see that
\begin{align*}
\widehat{\rho } _{b}(f, g) = \inf \{\varepsilon > 0 : \ \text{dist}(f(x),
g(B_{x}(\varepsilon ))) < \varepsilon \ \text{and} \ \text{dist}(g(x),
f(B_{x}(\varepsilon ))) < \varepsilon  \\
\text{for all} \ x \in X \}.
\end{align*}

(2) \ For any $ f, g \in
\mathfrak{M}(X \rightarrow Y), $ we define
\begin{align*}
\widehat{\rho } _{u}(f, g) = \inf \{\varepsilon > 0 : \
\text{there exists a} \ \delta > 0 \ \text{such that} \
\sup _{x^{\prime } \in B_{x}(\delta )} \widehat{d}_{Y} (f(x), g(x^{\prime }))
\leq \varepsilon \\
\text{and} \ \sup _{x^{\prime } \in B_{x}(\delta )}
\widehat{d}_{Y} (g(x), f(x^{\prime })) \leq \varepsilon
\ \text{for all} \ x \in X \} \ \text{if} \ f \neq g,
\end{align*}
$ \widehat{\rho } _{u}(f, g) = 0 $ if $ f = g. $ \\
It is easy to see that for $ f \neq g, $
\begin{align*}
\widehat{\rho } _{u}(f, g) = \inf \{\varepsilon > 0 : \
\text{there exists a} \ \delta > 0 \ \text{such that} \
\widehat{d}_{Y} (f(x), g(x^{\prime })) \leq \varepsilon  \\
\text{for all} \ x, x^{\prime } \in X \ \text{with} \
\widehat{d}_{X} (x, x^{\prime }) < \delta \}.
\end{align*}
\par     \vskip  0.2 cm

For $ f \in \mathfrak{M}(X \rightarrow Y) $ as above, recall that $ f $
is called $ Lipschitz $ if there exists a real number $ c \geq 0 $ such that
$ \widehat{d}_{Y} (f(x_{1}), f(x_{2})) \leq c \cdot
\widehat{d}_{X} (x_{1}, x_{2}) $ for all $ x_{1}, x_{2} \in X. $ The
$ dilatation $ of a Lipschitz map $ f $ is defined by
$$ \text{dil} f = \sup _{x_{1}, x_{2} \in X, \ x_{1} \neq x_{2}}
\frac{\widehat{d}_{Y} (f(x_{1}), f(x_{2}))}{\widehat{d}_{X} (x_{1}, x_{2})}
\quad \text{(see [BBI, p.9] and [G, p.1])}. $$
We denote $ \text{dil} f = + \infty $ if $ f $ is not a $ Lipschitz $ map.
\par     \vskip  0.2 cm

{\bf Theorem 3.5.} \ Let $ (X, \widehat{d}_{X} ) $ and $ (Y, \widehat{d}_{Y} ) $
be two non-Archimedean metric spaces. \\
(1) \ $ \widehat{\rho } _{u} $ is a non-Archimedean metric in
$ \mathfrak{M}(X \rightarrow Y), $ that is,
$ ( \mathfrak{M}(X \rightarrow Y), \ \widehat{\rho } _{u}) $ is
a non-Archimedean metric space.  \\
(2) \ For any $ f, g \in \mathfrak{M}(X \rightarrow Y), $ we have
$ \widehat{\rho } _{H}(f, g) \leq \widehat{\rho } _{u}(f, g). $ In particular,
if both $ f $ and $ g $ are nonexpanding, then
$ \widehat{\rho } _{H}(f, g) = \widehat{\rho } _{u}(f, g). $ \\
(3) \ For any $ f, g \in \mathfrak{M}(X \rightarrow Y), $
we have $ \widehat{\rho } _{H}(f, g) =
\widehat{\rho } _{b}(f, g). $ Particularly, if either $ f $ or $ g $ is
a Lipschitz map, and $ \min \{\text{dil} f, \ \text{dil} g \} \geq 1, $ then
$$ \widehat{\rho } _{s}(f, g) \leq \min \{\text{dil} f, \ \text{dil} g \} \cdot
\widehat{\rho } _{H}(f, g). $$
{\bf Proof.} \ (1) \ Obviously, $ \widehat{\rho } _{u} (f, g) =
\widehat{\rho } _{u} (g, f) \geq 0 $ for all
$ f, g \in \mathfrak{M}(X \rightarrow Y). $ If $ f \neq g, $ then
$ f(x_{0}) \neq g (x_{0}) $ for some $ x_{0} \in X, $ so
$ \widehat{d}_{Y} (f(x_{0}), g(x_{0})) > 0. $ Also for any $ \delta > 0,
\ \sup _{x^{\prime } \in B_{x_{0}}(\delta )} \widehat{d}_{Y} (f(x_{0}), g(x^{\prime }))
\geq \widehat{d}_{Y} (f(x_{0}), g(x_{0})). $ Hence, by definition,
$ \widehat{\rho } _{u} (f, g) \geq \widehat{d}_{Y} (f(x_{0}), g(x_{0})) > 0. $
Now we come to verify the strong triangle inequality.
To see this, let $ f, g, h \in \mathfrak{M}(X \rightarrow Y), $ we may
as well assume that $ \widehat{\rho } _{u} (f, h) \geq \widehat{\rho } _{u} (g, h), $
and we want to show that $ \widehat{\rho } _{u} (f, g) \leq \widehat{\rho } _{u} (f, h). $
In fact, for any $ \varepsilon > \widehat{\rho } _{u} (f, h), $ note that also
$ \varepsilon > \widehat{\rho } _{u} (g, h), $ it then follows easily by definition
that, there is a $ \delta > 0 $ such that
\begin{align*}
&\sup _{x^{\prime } \in B_{x}(\delta )} \widehat{d}_{Y} (f(x), h(x^{\prime }))
< \varepsilon , \quad
\sup _{x^{\prime } \in B_{x}(\delta )} \widehat{d}_{Y} (h(x), f(x^{\prime }))
< \varepsilon ,  \\
& \sup _{x^{\prime } \in B_{x}(\delta )} \widehat{d}_{Y} (g(x), h(x^{\prime }))
< \varepsilon , \quad  \sup _{x^{\prime } \in B_{x}(\delta )}
\widehat{d}_{Y} (h(x), g(x^{\prime })) < \varepsilon \
(\forall x \in X).
\end{align*}
Particularly, $ \widehat{d}_{Y} (g(x), h(x)) < \varepsilon \ (\forall x \in X). $
So, $$ \widehat{d}_{Y} (f(x), g(x^{\prime })) \leq \max
\{\widehat{d}_{Y} (f(x), h(x^{\prime })), \
\widehat{d}_{Y} (g(x^{\prime}), h(x^{\prime })) \} < \varepsilon \
(\forall x^{\prime } \in B_{x}(\delta )), $$ hence
$ \sup _{x^{\prime } \in B_{x}(\delta )} \widehat{d}_{Y} (f(x), g(x^{\prime }))
\leq \varepsilon , $ similarly,
$ \sup _{x^{\prime } \in B_{x}(\delta )} \widehat{d}_{Y} (g(x), f(x^{\prime }))
\leq \varepsilon \ (\forall x \in X), $ which implies
$ \widehat{\rho } _{u} (f, g) \leq \varepsilon . $ Since
$ \varepsilon > \widehat{\rho } _{u} (f, h) $ is arbitrary, we get
$ \widehat{\rho } _{u} (f, g) \leq \widehat{\rho } _{u} (f, h). $ So the strong
triangle inequality holds. This proves (1). \\
(2) \ For any $ \varepsilon > \widehat{\rho } _{u} (f, g), $ by
definition, it is easy to see that there is a positive number
$ \delta \leq \varepsilon $ such that $ \sup _{x^{\prime } \in B_{x}(\delta )}
\widehat{d}_{Y} (f(x), g(x^{\prime })) < \varepsilon \ (\forall x \in X). $
Then
\begin{align*}
\text{dist} ((x, f(x)), \Gamma _{g}) &= \inf _{x^{\prime \prime } \in X}
\max \{ \widehat{d}_{X} (x, x^{\prime \prime }), \
\widehat{d}_{Y} (f(x), g(x^{\prime \prime })) \}  \\
&\leq \sup _{x^{\prime } \in B_{x}(\delta )}
\max \{ \widehat{d}_{X} (x, x^{\prime }), \
\widehat{d}_{Y} (f(x), g(x^{\prime })) \} \\
&\leq \varepsilon \ (\forall x \in X).
\end{align*}
Similarly, $ \text{dist} ((x, g(x)), \Gamma _{f}) \leq \varepsilon \
(\forall x \in X). $ Therefore, by definition, it follows easily that
$ \widehat{\rho } _{H}(f, g) =
\widehat{d}_{Z, H} (\Gamma _{f}, \Gamma _{G}) \leq \varepsilon , $ which
implies $ \widehat{\rho } _{H}(f, g) \leq \widehat{\rho } _{u}(f, g) $
because $ \varepsilon > \widehat{\rho } _{u} (f, g) $ is arbitrary.
Now we assume that $ f $ and $ g $ are nonexpanding maps. Let
$ \varepsilon > \widehat{\rho } _{H}(f, g), $ by definition,
$ \varepsilon > \text{dist} ((x, f(x)), \Gamma _{g}) $ for all $ x \in X, $
which implies $ \widehat{d}_{Y} (f(x), g(x^{\prime })) < \varepsilon $
for some $ x^{\prime} \in B_{x}(\varepsilon ) \ (\forall x \in X). $
Since $ g $ is non-expanding, for any $ x^{\prime \prime } \in B_{x}(\varepsilon ), $
we have
\begin{align*}
\widehat{d}_{Y} (f(x), g(x^{\prime \prime })) &\leq
\max \{ \widehat{d}_{Y} (f(x), g(x^{\prime })), \
\widehat{d}_{Y} (g(x^{\prime }), g(x^{\prime \prime })) \} \\
&\leq \max \{ \widehat{d}_{Y} (f(x), g(x^{\prime })), \
\widehat{d}_{X} (x^{\prime }, x^{\prime \prime }) \} < \varepsilon .
\end{align*}
Hence $ \sup _{x^{\prime \prime } \in B_{x}(\varepsilon )}
\widehat{d}_{Y} (f(x), g(x^{\prime \prime })) \leq \varepsilon , $ similarly,
$ \sup _{x^{\prime \prime } \in B_{x}(\varepsilon )}
\widehat{d}_{Y} (g(x), f(x^{\prime \prime })) \leq \varepsilon \
(\forall x \in X). $ So by definition, $ \widehat{\rho } _{u}(f, g)
\leq \varepsilon , $ which implies $ \widehat{\rho } _{u}(f, g) \leq
\widehat{\rho } _{H}(f, g) $ because $ \varepsilon > \widehat{\rho } _{H}(f, g) $
is arbitrary. Therefore, $ \widehat{\rho } _{u}(f, g) =
\widehat{\rho } _{H}(f, g), $ this proves (2). \\
(3) \ Firstly, we verify the equality. For any
$ \varepsilon > \widehat{\rho } _{b}(f, g), $ by definition, there exists a
positive number $ \eta < \varepsilon $ such that, for every $ x \in X,
\ \widehat{d}_{Y} (f(x), g(x^{\prime }))
< \eta $ and $ \widehat{d}_{Y} (g(x), f(x^{\prime \prime }))
< \eta $ for some $ x^{\prime }, x^{\prime \prime }
\in B_{x}(\eta ), $ which implies that
\begin{align*}
&\text{dist}((x, f(x)), \Gamma _{g}) = \inf _{x_{2} \in X}
\max \{ \widehat{d}_{X} (x, x_{2}), \
\widehat{d}_{Y} (f(x), g(x_{2})) \} < \eta < \varepsilon
\quad \text{and} \\
&\text{dist}((x, g(x)), \Gamma _{f}) = \inf _{x_{1} \in X}
\max \{ \widehat{d}_{X} (x_{1}, x), \
\widehat{d}_{Y} (f(x_{1}), g(x)) \} < \eta < \varepsilon .
\end{align*}
Hence by definition, $ \widehat{\rho } _{H}(f, g) =
\widehat{d}_{X \times Y, H} (\Gamma _{f}, \Gamma _{g}) < \varepsilon $
for all $ \varepsilon > \widehat{\rho } _{b}(f, g), $ and so
$ \widehat{\rho } _{H}(f, g) \leq \widehat{\rho } _{b}(f, g). $ \ Conversely, let
$ \varepsilon > \widehat{\rho } _{H}(f, g), $ then, for every $ x \in X, $
there exist $ x^{\prime }, x^{\prime \prime } \in X $ such that
$$ \max \{ \widehat{d}_{X} (x, x^{\prime }), \
\widehat{d}_{Y} (f(x), g(x^{\prime })) \} < \varepsilon \ \text{and} \
\max \{ \widehat{d}_{X} (x, x^{\prime \prime }), \
\widehat{d}_{Y} (g(x), f(x^{\prime \prime })) \} < \varepsilon , $$
so $ x^{\prime }, x^{\prime \prime } \in B_{x}(\varepsilon ), \
\widehat{d}_{Y} (f(x), g(x^{\prime })) < \varepsilon $ and
$ \widehat{d}_{Y} (g(x), f(x^{\prime \prime })) < \varepsilon . $
Then by definition, $ \varepsilon \geq \widehat{\rho } _{b}(f, g), $
and so $ \widehat{\rho } _{b}(f, g) \leq \widehat{\rho } _{H}(f, g) $ because
$ \varepsilon > \widehat{\rho } _{H}(f, g) $ is arbitrary. Therefore,
$ \widehat{\rho } _{b}(f, g) = \widehat{\rho } _{H}(f, g). $ \\
Next, we come to verify the inequality. We may assume that
$ 1 \leq \text{dil} f \leq \text{dil} g. $ Denote $ c = \text{dil} f, $
we only need to show that $ \widehat{\rho } _{s}(f, g) \leq
c \cdot \widehat{\rho } _{b}(f, g). $
To see this, let $ \varepsilon > \widehat{\rho } _{b}(f, g), $ by definition, for
every $ x \in X, $ there exists a $ x^{\prime } \in B_{x}(\varepsilon ) $
such that $ \widehat{d}_{Y} (g(x), f(x^{\prime }))
< \varepsilon . $ Then by the definition of $ \text{dil} f, $ we get
$$ \widehat{d}_{Y} (g(x), f(x)) \leq \max \{\widehat{d}_{Y} (g(x), f(x^{\prime })),
\ \widehat{d}_{Y} (f(x^{\prime }), f(x)) \}
< \max \{\varepsilon , \ c \cdot \varepsilon \}
= c \cdot \varepsilon , $$ which implies that $ \widehat{\rho } _{s}(f, g) \leq
c \cdot \varepsilon , $ and so $ \widehat{\rho } _{s}(f, g) \leq c \cdot
\widehat{\rho } _{b}(f, g) $ because
$ \varepsilon > \widehat{\rho } _{b}(f, g) $ is arbitrary.
This proves (3), and the proof of Theorem 3.5
is completed. \quad $ \Box $
\par   \vskip  0.2 cm
Note that, for $ f, g $ as in the above Theorem 3.5(3), if
$ \min \{\text{dil} f, \ \text{dil} g \} \leq 1, $ we may as well assume
that $ \text{dil} f \leq 1, $ then $ f $ is a non-expanding map, so
from the above Remark 3.3, we know that $ \widehat{\rho } _{s}(f, g) =
\widehat{\rho } _{H}(f, g). $

\par     \vskip  0.3 cm

\hspace{-0.6cm}{\bf 4. \ Ultrametric structures on ball-type
mappings}

\par  \vskip  0.2 cm

Let $ (X, \widehat{d}_{X} ) $ and $ (Y, \widehat{d}_{Y} ) $
be two non-Archimedean metric spaces, recall that
$ \mathcal{M}_{\flat } (X) = \{ \text{ all balls in } X \} $
is a non-Archimedean metric space with the metric $ \widehat{d}_{X, H} $
as before. We denote $$ \mathfrak{D}_{\flat }(X, Y) =
\{\text{all maps} \ P : \mathcal{M}_{\flat } (X)
\rightarrow Y \} = \mathfrak{M}(\mathcal{M}_{\flat } (X) \rightarrow Y). $$
for any $ P_{1}, P_{2} \in \mathfrak{D}_{\flat }(X, Y), $ recall that
$$ \widehat{\rho} _{H}(P_{1}, P_{2}) = \widehat{d}_{W, H} (\Gamma _{P_{1}},
\Gamma _{P_{2}}), \ \ \widehat{\rho } _{s}(P_{1}, P_{2}) = \sup _{B
\in \mathcal{M}_{\flat } (X)} \widehat{d}_{Y} (P_{1}(B), P_{2}(B)),
$$ where $ W = \mathcal{M}_{\flat } (X) \times Y $ is endowed the
square metric $ \widehat{d}_{W} $ of $ \widehat{d}_{X, H} $ and $
\widehat{d}_{Y} $ defined before. Let $ \varepsilon > 0, $ for
convenience, we will write $ B^{\varepsilon } = U_{\varepsilon} (B),
$ the $ \varepsilon -$neighborhood of a ball $ B $ in $ X. $ Note
that by Lemma 2.3 above, $ B^{\varepsilon } $ is also a ball in $ X.
$ For all spaces considered in the following, we assume that every
ball in them contains at least two distinct elements.
\par     \vskip  0.2 cm

{\bf Definition 4.1.} \ Let $ (X, \widehat{d}_{X} ) $ and $ (Y, \widehat{d}_{Y} ) $
be two non-Archimedean metric spaces. For any $ P_{1}, P_{2}
\in \mathfrak{D}_{\flat }(X, Y) $ and $ \lambda > 0, $ we define
\begin{align*}
\widehat{\beta }_{X, Y}^{\lambda }(P_{1}, P_{2}) =
\inf \{\varepsilon > 0 : \ \widehat{d}_{Y} (P_{1}(B), P_{2}(B^{\lambda \varepsilon }))
\leq \varepsilon \ \text{and} \ \widehat{d}_{Y} (P_{2}(B), P_{1}(B^{\lambda \varepsilon }))
\leq \varepsilon \ \\
\text{for all} \ B \in \mathcal{M}_{\flat } (X) \}
\ \text{if} \ P_{1} \neq P_{2},
\end{align*}
$ \widehat{\beta }_{X, Y}^{\lambda }(P_{1}, P_{2}) = 0 $
if $ P_{1} = P_{2}; $
\begin{align*}
\widehat{\beta }_{X, Y}^{\ast \lambda }(P_{1}, P_{2})=
\inf \{\varepsilon > 0 : \ \text{for every} \ B \in \mathcal{M}_{\flat } (X),
\ \text{there exist positive numbers} \\
\varepsilon _{B}, \varepsilon _{B}^{\prime } \leq \varepsilon \
\text{such that} \
\widehat{d}_{Y} (P_{1}(B), P_{2}(B^{\lambda \varepsilon _{B}})) \leq \varepsilon
\ \text{and} \ \widehat{d}_{Y} (P_{2}(B), P_{1}(B^{\lambda \varepsilon _{B}^{\prime }}))
\leq \varepsilon \}.
\end{align*}
For the case $ \lambda = 1, $ we write $ \widehat{\beta }_{X, Y}^{1}
= \widehat{\beta }_{X, Y} $ and $ \widehat{\beta }_{X, Y}^{\ast 1} =
\widehat{\beta }_{X, Y}^{\ast }. $ By definition, it is easy to see that
$ \widehat{\beta }_{X, Y}^{\ast \lambda } (P_{1}, P_{2})
\leq \widehat{\beta }_{X, Y}^{\lambda }(P_{1}, P_{2}) \
(\forall P_{1}, P_{2} \in \mathfrak{D}_{\flat }(X, Y)). $ Moreover,
for any $ P \in \mathfrak{D}_{\flat }(X, Y), $ if we
denote $$ \eta _{\lambda , P} = \inf \{\varepsilon > 0 : \
\widehat{d}_{Y} (P(B), P(B^{\lambda \varepsilon })) \leq \varepsilon
\ \text{for all} \ B \in \mathcal{M}_{\flat } (X) \}, $$ then it can be
easily verified that
$ \widehat{\beta }_{X, Y}^{\lambda }(P_{1}, P_{2}) \geq \max
\{ \eta _{\lambda , P_{1}}, \ \eta _{\lambda , P_{2}} \} $
for all $ P_{1}, P_{2} \in \mathfrak{D}_{\flat }(X, Y). $
We also denote
\begin{align*}
\mathfrak{D}_{\flat }^{(\lambda )}(X, Y) =
\{P \in \mathfrak{D}_{\flat }(X, Y) : \
\widehat{d}_{Y} (P(B), P(B^{\lambda \varepsilon })) \leq \varepsilon \\
\text{for all} \ \varepsilon > 0 \ \text{and all} \
B \in \mathcal{M}_{\flat } (X) \}.
\end{align*}
Obviously, $ \eta _{\lambda , P} = 0 $
for all $ P \in \mathfrak{D}_{\flat }^{(\lambda )}(X, Y). $
Note that, if $ P \in \mathfrak{D}_{\flat }(X, Y) $ is nonexpanding, then
$ \widehat{d}_{Y} (P(B), P(B^{\varepsilon })) \leq
\widehat{d}_{X, H}(B, B^{\varepsilon }) \leq \varepsilon , $ so
$ P \in \mathfrak{D}_{\flat }^{(1)}(X, Y). $ Hence
$ \mathfrak{D}_{\flat }^{(1)}(X, Y) $ contains all nonexpanding maps from
$ \mathcal{M}_{\flat } (X) $ to $ Y. $
\par   \vskip  0.2 cm

{\bf Theorem 4.2.} (1) \ The above defined functions
$ \widehat{\beta }_{X, Y}^{\lambda } $
and $ \widehat{\beta }_{X, Y}^{\ast \lambda } $ are non-Archimedean
metrics on $ \mathfrak{D}_{\flat }(X, Y). $
\par     \vskip  0.1 cm
(2) \ For any $ \lambda > 0 $ and $ P_{1}, P_{2}
\in \mathfrak{D}_{\flat }(X, Y), $ we have
$$ \widehat{\beta }_{X, Y}^{\ast \lambda }(P_{1}, P_{2}) \leq
\widehat{\rho } _{s}(P_{1}, P_{2}) \leq
\widehat{\beta }_{X, Y}^{\lambda }(P_{1}, P_{2}), \quad
\widehat{\rho} _{H}(P_{1}, P_{2}) \leq \max \{1, \lambda \} \cdot
\widehat{\beta }_{X, Y}^{\ast \lambda }(P_{1}, P_{2}). $$ In particular,
if $ P_{1}, P_{2} \in \mathfrak{D}_{\flat }^{(1)}(X, Y), $ then
$$ \widehat{\rho} _{H}(P_{1}, P_{2}) = \widehat{\rho } _{s}(P_{1}, P_{2})
= \widehat{\beta }_{X, Y}(P_{1}, P_{2}) =
\widehat{\beta }_{X, Y}^{\ast }(P_{1}, P_{2}). $$
(3) \ For $ P \in \mathfrak{D}_{\flat }^{(1)}(X, Y) $
and a sequence $ \{P_{n} \}_{n = 1}^{\infty } $ in $ \mathfrak{D}_{\flat }(X, Y), $
the following statements are equivalent: \\
(A) \ $ \widehat{\beta }_{X, Y}(P_{n}, P) \rightarrow 0 $ as $ n \rightarrow
\infty . $ \\
(B) \ $ \sup _{B \in \mathcal{M}_{\flat }(X)} \widehat{d}_{Y} (P_{n}(B), P(B))
\rightarrow 0 $ as $ n \rightarrow \infty , $ in other words, $ P_{n} $ converges
uniformly to $ P $ as $ n \rightarrow \infty . $
\par     \vskip  0.1 cm
{\bf Proof.} \ (1) \ We only verify $ \widehat{\beta }_{X, Y}^{\lambda }, $
the other is similar. By definition,
$ \widehat{\beta }_{X, Y}^{\lambda }(P_{1}, P_{2}) =
\widehat{\beta }_{X, Y}^{\lambda }(P_{2}, P_{1}) \geq 0 $ and
$ \widehat{\beta }_{X, Y}^{\lambda }(P_{1}, P_{1}) = 0 \ (\forall
P_{1}, P_{2} \in \mathfrak{D}_{\flat }(X, Y)). $ Now suppose $ P_{1} \neq P_{2} $
and $ \widehat{\beta }_{X, Y}^{\lambda }(P_{1}, P_{2}) = 0 $ for some
$ P_{1}, P_{2}, $ then, for any $ B \in \mathcal{M}_{\flat } (X) $ and
positive number $ \varepsilon \leq \text{diam}(B) / \lambda , $
it follows easily from the definition of
$ \widehat{\beta }_{X, Y}^{\lambda } $ that
$ \widehat{d}_{Y} (P_{1}(B), P_{2}(B^{\lambda \varepsilon ^{\prime }}))
\leq \varepsilon ^{\prime } $ for some positive number
$ \varepsilon ^{\prime } < \varepsilon . $ Since
$ \lambda \varepsilon ^{\prime } < \lambda \varepsilon \leq \text{diam}(B), $
by Lemma 2.3 above, we have
$ B^{\lambda \varepsilon ^{\prime }} = B, $ which implies
$ \widehat{d}_{Y} (P_{1}(B), P_{2}(B)) < \varepsilon , $ and then
$ \widehat{d}_{Y} (P_{1}(B), P_{2}(B)) = 0 $ because
$ \varepsilon \leq \text{diam}(B) / \lambda $ is arbitrary. Hence,
$ P_{1}(B) = P_{2}(B) \ (\forall B \in \mathcal{M}_{\flat } (X)), $ and so
$ P_{1} = P_{2}, $ a contradiction ! Therefore,
$ \widehat{\beta }_{X, Y}^{\lambda }(P_{1}, P_{2}) > 0 $ if $ P_{1} \neq P_{2}. $
Lastly, we come to verify the strong triangle inequality. Let $ P_{1},
P_{2}, P_{3} \in \mathfrak{D}_{\flat }(X, Y), $ for any $ \varepsilon >
\max \{\widehat{\beta }_{X, Y}^{\lambda }(P_{1}, P_{2}), \
\widehat{\beta }_{X, Y}^{\lambda }(P_{2}, P_{3}) \}, $ by definition,
there exist positive numbers $ \varepsilon _{1}, \varepsilon _{2}
< \varepsilon $ such that
\begin{align*}
&\widehat{d}_{Y} (P_{1}(B), P_{2}(B^{\lambda \varepsilon _{1}}))
\leq \varepsilon _{1}, \ \widehat{d}_{Y} (P_{2}(B), P_{1}(B^{\lambda \varepsilon _{1}}))
\leq \varepsilon _{1}; \\
&\widehat{d}_{Y} (P_{2}(B), P_{3}(B^{\lambda \varepsilon _{2}}))
\leq \varepsilon _{2}, \ \widehat{d}_{Y} (P_{3}(B), P_{2}(B^{\lambda \varepsilon _{2}}))
\leq \varepsilon _{2} \ (\forall B \in \mathcal{M}_{\flat }(X)).
\end{align*}
We may as well assume that $ \varepsilon _{1} \leq \varepsilon _{2}. $ Then
by Lemma 2.3 above, $ (B^{\lambda \varepsilon _{2}})^{\lambda \varepsilon _{1}}
= (B^{\lambda \varepsilon _{1}})^{\lambda \varepsilon _{2}} =
B^{\lambda \varepsilon _{2}}, $ so
\begin{align*}
&\widehat{d}_{Y} (P_{2}(B^{\lambda \varepsilon _{2}}), P_{1}(B^{\lambda \varepsilon _{2}}))
= \widehat{d}_{Y} (P_{2}(B^{\lambda \varepsilon _{2}}),
P_{1}((B^{\lambda \varepsilon _{2}})^{\lambda \varepsilon _{1}}))
\leq \varepsilon _{1}; \\
&\widehat{d}_{Y} (P_{2}(B^{\lambda \varepsilon _{1}}), P_{3}(B^{\lambda \varepsilon _{2}}))
= \widehat{d}_{Y} (P_{2}(B^{\lambda \varepsilon _{1}}),
P_{3}((B^{\lambda \varepsilon _{1}})^{\lambda \varepsilon _{2}}))
\leq \varepsilon _{2}.
\end{align*}
Thus $ \widehat{d}_{Y} (P_{1}(B), P_{3}(B^{\lambda \varepsilon _{2}})) \leq
\max \{\widehat{d}_{Y} (P_{1}(B), P_{2}(B^{\lambda \varepsilon _{1}})),
\widehat{d}_{Y} (P_{2}(B^{\lambda \varepsilon _{1}}),
P_{3}(B^{\lambda \varepsilon _{2}})) \} \leq \max \{\varepsilon _{1},
\varepsilon _{2} \} \leq \varepsilon _{2}. $ Similarly,
$ \widehat{d}_{Y} (P_{3}(B), P_{1}(B^{\lambda \varepsilon _{2}})) \leq
\varepsilon _{2}. $ Therefore, by definition, we obtain that
$ \varepsilon > \varepsilon _{2} \geq
\widehat{\beta }_{X, Y}^{\lambda }(P_{1}, P_{3}), $ and so \\
$ \widehat{\beta }_{X, Y}^{\lambda }(P_{1}, P_{3}) \leq
\max \{\widehat{\beta }_{X, Y}^{\lambda }(P_{1}, P_{2}), \
\widehat{\beta }_{X, Y}^{\lambda }(P_{2}, P_{3}) \} $ by the
arbitrary choice of $ \varepsilon . $ This proves (1).
\par     \vskip  0.1 cm
(2) \ We first verify $ \widehat{\beta }_{X, Y}^{\ast \lambda }(P_{1}, P_{2})
\leq \widehat{\rho } _{s}(P_{1}, P_{2}). $ We may assume $ P_{1} \neq P_{2}, $
and write $ \varepsilon _{0} = \widehat{\rho } _{s}(P_{1}, P_{2}), $ then
$ \varepsilon _{0} > 0. $ For every $ B \in \mathcal{M}_{\flat }(X), $ let
$ \varepsilon _{B} = \min \{ \text{diam} B / \lambda , \varepsilon _{0} \}, $
then by Lemma 2.3, we have $ B^{\lambda \varepsilon _{B}} = B. $ So
\begin{align*}
&\widehat{d}_{Y} (P_{1}(B), P_{2}(B^{\lambda \varepsilon _{B}})) =
\widehat{d}_{Y} (P_{1}(B), P_{2}(B)) \leq \widehat{\rho } _{s}(P_{1}, P_{2})
 = \varepsilon _{0}, \\
&\widehat{d}_{Y} (P_{2}(B), P_{1}(B^{\lambda \varepsilon _{B}})) =
\widehat{d}_{Y} (P_{2}(B), P_{1}(B)) \leq \widehat{\rho } _{s}(P_{1}, P_{2})
 = \varepsilon _{0},
\end{align*}
and the inequality follows by definition. \\
Next we verify $ \widehat{\rho } _{s}(P_{1}, P_{2}) \leq
\widehat{\beta }_{X, Y}^{\lambda }(P_{1}, P_{2}). $ We may assume $ P_{1} \neq P_{2}. $
For any $ \varepsilon > \widehat{\beta }_{X, Y}^{\lambda }(P_{1}, P_{2}), $
by definition, there exists a positive number $ \varepsilon ^{\prime } <
\varepsilon $ such that $$ \widehat{d}_{Y} (P_{1}(B),
P_{2}(B^{\lambda \varepsilon ^{\prime }})) \leq \varepsilon ^{\prime }
\ \text{and} \ \widehat{d}_{Y} (P_{2}(B),
P_{1}(B^{\lambda \varepsilon ^{\prime }})) \leq \varepsilon ^{\prime }
\ (\forall B \in \mathcal{M}_{\flat }(X)). $$ In particular,
$ \widehat{d}_{Y} (P_{2}(B^{\lambda \varepsilon ^{\prime }}),
P_{1}(B^{\lambda \varepsilon ^{\prime }})) =
\widehat{d}_{Y} (P_{2}(B^{\lambda \varepsilon ^{\prime }}),
P_{1}((B^{\lambda \varepsilon ^{\prime }})^{\lambda \varepsilon ^{\prime }}))
\leq \varepsilon ^{\prime } $ (see Lemma 2.3 above). So $$ \widehat{d}_{Y} (P_{1}(B),
P_{1}(B^{\lambda \varepsilon ^{\prime }})) \leq \max \{\widehat{d}_{Y} (P_{1}(B),
P_{2}(B^{\lambda \varepsilon ^{\prime }})), \
\widehat{d}_{Y} (P_{2}(B^{\lambda \varepsilon ^{\prime }}),
P_{1}(B^{\lambda \varepsilon ^{\prime }})) \} \leq \varepsilon ^{\prime }, $$
and then $ \widehat{d}_{Y} (P_{1}(B), P_{2}(B)) \leq \max \{\widehat{d}_{Y} (P_{1}(B),
P_{1}(B^{\lambda \varepsilon ^{\prime }})), \
\widehat{d}_{Y} (P_{1}(B^{\lambda \varepsilon ^{\prime }}),
P_{2}(B)) \} \leq \varepsilon ^{\prime } < \varepsilon \
(\forall B \in \mathcal{M}_{\flat }(X)), $ which implies
$ \widehat{\rho } _{s}(P_{1}, P_{2}) \leq \varepsilon ^{\prime } < \varepsilon , $
and so $ \widehat{\rho } _{s}(P_{1}, P_{2}) \leq
\widehat{\beta }_{X, Y}^{\lambda }(P_{1}, P_{2}) $ by the arbitrary choice of
$ \varepsilon . $ \\
Next we verify $ \widehat{\rho} _{H}(P_{1}, P_{2}) \leq \max \{1, \lambda \}
\cdot \widehat{\beta }_{X, Y}^{\ast \lambda }(P_{1}, P_{2}). $ For any
$ \varepsilon > \widehat{\beta }_{X, Y}^{\ast \lambda }(P_{1}, P_{2}), $
by definition, there is a positive number $ \varepsilon _{1} < \varepsilon $
such that $ \widehat{d}_{Y} (P_{1}(B), P_{2}(B^{\lambda \varepsilon _{B}})
\} \leq \varepsilon _{1} $ and $ \widehat{d}_{Y} (P_{2}(B),
P_{1}(B^{\lambda \varepsilon _{B}^{\prime }}) \} \leq \varepsilon _{1}
$ for some positive numbers $ \varepsilon _{B}, \varepsilon _{B}^{\prime }
\leq \varepsilon _{1} \ (\forall B \in \mathcal{M}_{\flat }(X)). $
By Lemma 2.1 and Lemma 2.3 above, we have $ \widehat{d}_{X, H}
(B, B^{\lambda \varepsilon _{B}}) \leq \lambda \varepsilon _{B}
< \lambda \varepsilon $ and $ \widehat{d}_{X, H}
(B, B^{\lambda \varepsilon _{B}^{\prime }}) \leq \lambda
\varepsilon _{B}^{\prime } < \lambda \varepsilon , $ so
$ B^{\lambda \varepsilon _{B}}, B^{\lambda \varepsilon _{B}^{\prime }}
\in B_{B}(\lambda \varepsilon ), $ a ball in the space
$ \mathcal{M}_{\flat }(X). $ Hence by Definition 3.4(1) above,
$ \widehat{\rho} _{b} (P_{1}, P_{2}) \leq \max \{1, \lambda \} \varepsilon . $
Since $ \varepsilon $ is arbitrary, by Theorem 3.5(3) above, we obtain
$ \widehat{\rho} _{H}(P_{1}, P_{2}) = \widehat{\rho} _{b} (P_{1}, P_{2})
\leq \max \{1, \lambda \}
\cdot \widehat{\beta }_{X, Y}^{\ast \lambda }(P_{1}, P_{2}). $ \\
Lastly, assume $ P_{1}, P_{2} \in \mathfrak{D}_{\flat }^{(1)}(X, Y), $
we come to verify the equalities. By the above discussion and Theorem 3.5(3),
we only need to show that
$ \widehat{\beta }_{X, Y}(P_{1}, P_{2}) \leq \widehat{\rho } _{b}(P_{1}, P_{2}). $
To see this, for any $ \varepsilon > \widehat{\rho} _{b}(P_{1}, P_{2}), $ by
definition, there is a positive number $ \varepsilon _{1} < \varepsilon $
such that, for every $ B \in \mathcal{M}_{\flat }(X), $ there exist
$ B_{1}, B_{2} \in \mathcal{M}_{\flat }(X) $ satisfying
$ \widehat{d}_{X, H} (B, B_{1}) < \varepsilon _{1}, \
\widehat{d}_{X, H} (B, B_{2}) < \varepsilon _{1}, \
\widehat{d}_{Y} (P_{1}(B), P_{2}(B_{1})) < \varepsilon _{1} $
and $ \widehat{d}_{Y} (P_{2}(B), P_{1}(B_{2}))
< \varepsilon _{1}. $ Then by the above Lemma 2.1 and Lemma 2.3, it can be easily
verified that $ B_{1}^{\varepsilon } = B_{2}^{\varepsilon } = B^{\varepsilon }. $
Since $ P_{1}, P_{2} \in \mathfrak{D}_{\flat }^{(1)}(X, Y), $ we have
\begin{align*}
&\widehat{d}_{Y} (P_{1}(B_{2}), P_{1}(B^{\varepsilon })) =
\widehat{d}_{Y} (P_{1}(B_{2}), P_{1}(B_{2}^{\varepsilon })) \leq \varepsilon , \\
&\widehat{d}_{Y} (P_{2}(B_{1}), P_{2}(B^{\varepsilon })) =
\widehat{d}_{Y} (P_{2}(B_{1}), P_{2}(B_{1}^{\varepsilon })) \leq \varepsilon .
\quad \text{So} \\
&\widehat{d}_{Y} (P_{1}(B), P_{2}(B^{\varepsilon })) \leq
\max \{\widehat{d}_{Y} (P_{1}(B), P_{2}(B_{1})), \
\widehat{d}_{Y} (P_{2}(B_{1}), P_{2}(B^{\varepsilon })) \} \leq \varepsilon
\ \text{and}  \\
&\widehat{d}_{Y} (P_{2}(B), P_{1}(B^{\varepsilon })) \leq
\max \{\widehat{d}_{Y} (P_{2}(B), P_{1}(B_{2})), \
\widehat{d}_{Y} (P_{1}(B_{2}), P_{1}(B^{\varepsilon })) \} \leq \varepsilon .
\end{align*}
Hence by definition, $ \widehat{\beta }_{X, Y}(P_{1}, P_{2}) \leq \varepsilon . $
Therefore $ \widehat{\beta }_{X, Y}(P_{1}, P_{2}) \leq
\widehat{\rho } _{b}(P_{1}, P_{2}) $ because
$ \varepsilon > \widehat{\rho} _{b}(P_{1}, P_{2}) $ is arbitrary.
This proves (2).
\par     \vskip  0.1 cm
(3) \ (A) $ \Rightarrow $ (B). For any $ \varepsilon > 0, $ there is
a positive integer $ N $ such that, for every integer $ n > N $ we have
$ \widehat{\beta }_{X, Y}(P_{n}, P) < \varepsilon . $ By definition, there exists
a positive number $ \varepsilon ^{\prime } < \varepsilon $ such that
$ \widehat{d}_{Y} (P_{n}(B), P(B^{\varepsilon ^{\prime }})) \leq
\varepsilon ^{\prime } \ (\forall B \in \mathcal{M}_{\flat }(X)). $ Since
$ P \in \mathfrak{D}_{\flat }^{(1)}(X, Y), $ we have
$ \widehat{d}_{Y} (P(B), P(B^{\varepsilon ^{\prime }})) \leq
\varepsilon ^{\prime }, $ so $$ \widehat{d}_{Y} (P_{n}(B), P(B)) \leq
\max \{\widehat{d}_{Y} (P_{n}(B), P(B^{\varepsilon ^{\prime }})), \
\widehat{d}_{Y} (P(B), P(B^{\varepsilon ^{\prime }})) \} \leq
\varepsilon ^{\prime }. $$ Therefore
$ \sup _{B \in \mathcal{M}_{\flat }(X)} \widehat{d}_{Y} (P_{n}(B), P(B))
\leq \varepsilon ^{\prime } < \varepsilon . $ \\
(B) $ \Rightarrow $ (A). For any $ \varepsilon > 0, $ there is
a positive integer $ N $ such that, for every integer $ n > N $ we have
$ \sup _{B \in \mathcal{M}_{\flat }(X)} \widehat{d}_{Y} (P_{n}(B), P(B))
< \varepsilon . $ Since $ P \in \mathfrak{D}_{\flat }^{(1)}(X, Y), $ we have
$ \widehat{d}_{Y} (P(B), P(B^{\varepsilon })) \leq \varepsilon \
(\forall B \in \mathcal{M}_{\flat }(X)). $ So
\begin{align*}
&\widehat{d}_{Y} (P_{n}(B), P(B^{\varepsilon })) \leq
\max \{\widehat{d}_{Y} (P_{n}(B), P(B)), \
\widehat{d}_{Y} (P(B), P(B^{\varepsilon })) \} \leq \varepsilon
\ \text{and}  \\
&\widehat{d}_{Y} (P(B), P_{n}(B^{\varepsilon })) \leq
\max \{\widehat{d}_{Y} (P(B), P(B^{\varepsilon })), \
\widehat{d}_{Y} (P(B^{\varepsilon }), P_{n}(B^{\varepsilon })) \}
\leq \varepsilon .
\end{align*}
Hence by definition, $ \widehat{\beta }_{X, Y}(P_{n}, P) \leq \varepsilon . $
This proves (3), and the proof of
Theorem 4.2 is completed. \quad $ \Box $
\par     \vskip  0.2 cm

Note that for $ P, P_{n} $ in Theorem 4.2.(3) above, if
$ P \in \mathfrak{D}_{\flat }(X, Y), $ then \\
$ \widehat{\beta }_{X, Y}(P_{n}, P) \rightarrow 0 \ \Rightarrow \
\widehat{d}_{Y} (P_{n}(B), P(B)) \rightarrow 0 \
(\forall B \in \mathcal{M}_{\flat }(X)). $
\par   \vskip  0.2 cm

{\bf Lemma 4.3.} \ On the interval $ (0, \infty ), $ the function $
\widehat{\beta }_{X, Y}^{\lambda } $ is increasing in the variable $ \lambda $
and $ \widehat{\beta }_{X, Y}^{\ast \lambda } $
is decreasing in this variable.
\par     \vskip  0.1 cm
{\bf Proof.} \ We only verify the case $ \widehat{\beta }_{X, Y}^{\lambda }, $
the other can be similarly done. Let $ \lambda _{1}, \lambda _{2} \in
(0, \infty ) $ with $ \lambda _{1} < \lambda _{2}. $ For any $ P_{1}, P_{2}
\in \mathfrak{D}_{\flat }(X, Y) $ with $ P_{1} \neq P_{2}, $ let
$ \varepsilon > \widehat{\beta }_{X, Y}^{\lambda _{2}} (P_{1}, P_{2}), $
by definition and Theorem 4.2.(2), there is a positive number
$ \varepsilon ^{\prime} < \varepsilon $ such that $ \varepsilon ^{\prime}
\geq \widehat{\rho } _{s}(P_{1}, P_{2}), \
\widehat{d}_{Y} (P_{1}(B), P_{2}(B^{\lambda_{2} \varepsilon ^{\prime}}))
\leq \varepsilon ^{\prime} $ and
$ \widehat{d}_{Y} (P_{2}(B), P_{1}(B^{\lambda_{2} \varepsilon ^{\prime}}))
\leq \varepsilon ^{\prime} \ (\forall B \in \mathcal{M}_{\flat }(X)). $
Note that by Lemma 2.3,
$ B^{\lambda_{1} \varepsilon ^{\prime}},
(B^{\lambda_{1} \varepsilon ^{\prime}})^{\lambda_{2} \varepsilon ^{\prime}}
= B^{\lambda_{2} \varepsilon ^{\prime}} \in \mathcal{M}_{\flat }(X), $
so
\begin{align*}
&\widehat{d}_{Y} (P_{1}(B^{\lambda_{1} \varepsilon ^{\prime}}),
P_{2}(B^{\lambda_{2} \varepsilon ^{\prime}})) =
\widehat{d}_{Y} (P_{1}(B^{\lambda_{1} \varepsilon ^{\prime}}),
P_{2}((B^{\lambda_{1} \varepsilon ^{\prime}})^{\lambda_{2} \varepsilon ^{\prime}}))
\leq \varepsilon ^{\prime}, \\
&\widehat{d}_{Y} (P_{2}(B^{\lambda_{1} \varepsilon ^{\prime}}),
P_{1}(B^{\lambda_{2} \varepsilon ^{\prime}})) =
\widehat{d}_{Y} (P_{2}(B^{\lambda_{1} \varepsilon ^{\prime}}),
P_{1}((B^{\lambda_{1} \varepsilon ^{\prime}})^{\lambda_{2} \varepsilon ^{\prime}}))
\leq \varepsilon ^{\prime}.
\end{align*}
Also by definition,
$ \widehat{d}_{Y} (P_{1}(B^{\lambda_{1} \varepsilon ^{\prime}}),
P_{2}(B^{\lambda_{1} \varepsilon ^{\prime}})) \leq
\widehat{\rho } _{s}(P_{1}, P_{2}) \leq \varepsilon ^{\prime}. $ Hence
\begin{align*}
&\widehat{d}_{Y} (P_{1}(B), P_{2}(B^{\lambda_{1} \varepsilon ^{\prime}}))
\leq \\
&\max
\{ \widehat{d}_{Y} (P_{1}(B), P_{2}(B^{\lambda_{2} \varepsilon ^{\prime}})),
\widehat{d}_{Y} (P_{2}(B^{\lambda_{2} \varepsilon ^{\prime}}),
P_{1}(B^{\lambda_{1} \varepsilon ^{\prime}})),
\widehat{d}_{Y} (P_{1}(B^{\lambda_{1} \varepsilon ^{\prime}}),
P_{2}(B^{\lambda_{1} \varepsilon ^{\prime}})) \} \leq \varepsilon ^{\prime}.
\end{align*}
Similarly,
$ \widehat{d}_{Y} (P_{2}(B), P_{1}(B^{\lambda_{1} \varepsilon ^{\prime}}))
\leq \varepsilon ^{\prime}. $ Therefore, by definition,
$ \widehat{\beta }_{X, Y}^{\lambda _{1}} (P_{1}, P_{2})
\leq \varepsilon ^{\prime} < \varepsilon , $ which implies
$ \widehat{\beta }_{X, Y}^{\lambda _{1}} (P_{1}, P_{2}) \leq
\widehat{\beta }_{X, Y}^{\lambda _{2}} (P_{1}, P_{2}) $ by the arbitrary
choice of $ \varepsilon . $ The proof is completed. \quad $ \Box $
\par    \vskip  0.2 cm

Now for any $ P_{1}, P_{2} \in \mathfrak{D}_{\flat }(X, Y), $ we denote
\begin{align*}
&\widehat{\beta }_{X, Y}^{0} (P_{1}, P_{2}) = \inf _{\lambda > 0}
\widehat{\beta }_{X, Y}^{\lambda } (P_{1}, P_{2}), \quad
\widehat{\beta }_{X, Y}^{\infty } (P_{1}, P_{2}) = \sup _{\lambda > 0}
\widehat{\beta }_{X, Y}^{\lambda } (P_{1}, P_{2}), \\
&\widehat{\beta }_{X, Y}^{\ast 0}(P_{1}, P_{2}) = \sup _{\lambda > 0}
\widehat{\beta }_{X, Y}^{\ast \lambda }(P_{1}, P_{2}), \quad
\widehat{\beta }_{X, Y}^{\ast \infty }(P_{1}, P_{2}) = \inf _{\lambda > 0}
\widehat{\beta }_{X, Y}^{\ast \lambda }(P_{1}, P_{2}).
\end{align*}
By Lemma 4.3 and Theorem 4.2 above, we have
\begin{align*}
&\widehat{\beta }_{X, Y}^{0} (P_{1}, P_{2}) =
\lim _{\lambda \rightarrow 0^{+}}
\widehat{\beta }_{X, Y}^{\lambda } (P_{1}, P_{2}), \quad
\widehat{\beta }_{X, Y}^{\infty } (P_{1}, P_{2}) =
\lim _{\lambda \rightarrow \infty }
\widehat{\beta }_{X, Y}^{\lambda } (P_{1}, P_{2}), \\
&\widehat{\beta }_{X, Y}^{\ast 0}(P_{1}, P_{2}) =
\lim _{\lambda \rightarrow 0^{+}}
\widehat{\beta }_{X, Y}^{\ast \lambda }(P_{1}, P_{2}), \quad
\widehat{\beta }_{X, Y}^{\ast \infty }(P_{1}, P_{2}) =
\lim _{\lambda \rightarrow \infty }
\widehat{\beta }_{X, Y}^{\ast \lambda }(P_{1}, P_{2}), \\
&\widehat{\beta }_{X, Y}^{\ast \infty }(P_{1}, P_{2}) \leq
\widehat{\beta }_{X, Y}^{\ast 0}(P_{1}, P_{2})
\leq \widehat{\rho } _{s}(P_{1}, P_{2}) \leq
\widehat{\beta }_{X, Y}^{0} (P_{1}, P_{2}) \leq
\widehat{\beta }_{X, Y}^{\infty } (P_{1}, P_{2}).
\end{align*}
By Theorem 4.2 above, it is not difficult to show that all
$ \widehat{\beta }_{X, Y}^{0}, \widehat{\beta }_{X, Y}^{\infty },
\widehat{\beta }_{X, Y}^{\ast 0} $ and
$ \widehat{\beta }_{X, Y}^{\ast \infty } $ are non-Archimedean
metrics on $ \mathfrak{D}_{\flat }(X, Y). $
\par    \vskip  0.2 cm

For the set $ \mathfrak{D}_{\flat }^{(\lambda )}(X, Y) $ defined in
Definition 4.1 above, it is easy to verify that \\
$ \mathfrak{D}_{\flat }^{(\lambda _{2})}(X, Y)
\subset \mathfrak{D}_{\flat }^{(\lambda _{1})}(X, Y) $ for any positive numbers
$ \lambda _{1} \leq \lambda _{2}. $ So particularly, we have a descending chain \
$ \mathfrak{D}_{\flat }(X, Y) \supset
\mathfrak{D}_{\flat }^{(1)}(X, Y) \supset
\mathfrak{D}_{\flat }^{(2)}(X, Y) \supset \cdots . $ \
We denote
\begin{align*}
&\mathfrak{D}_{\flat }^{(\infty )}(X, Y) = \bigcap _{ \lambda > 0}
\mathfrak{D}_{\flat }^{(\lambda )}(X, Y), \quad
\mathfrak{D}_{\flat }^{(0)}(X, Y) = \bigcup _{\lambda > 0}
\mathfrak{D}_{\flat }^{(\lambda )}(X, Y), \quad \text{and} \\
&\mathfrak{D}_{\flat }(X, Y)_{0} = \{ P \in \mathfrak{D}_{\flat }(X, Y) : \
P \ \text{is a constant map from} \ \mathcal{M}_{\flat } (X) \ \text{to} \
Y \}.
\end{align*}
Obviously, $ \mathfrak{D}_{\flat }(X, Y)_{0} \subset
\mathfrak{D}_{\flat }^{(\infty )}(X, Y). $
\par     \vskip  0.2 cm

{\bf Definition 4.4.} \ Let $ P, P_{1}, P_{2} \in
\mathfrak{D}_{\flat }(X, Y) $ and $ \varepsilon > 0. $ \\
(1) \ If there exists an $ \delta > 0 $ such that
$ \widehat{d}_{Y} (P(B), P(B^{\delta })) \leq \varepsilon $ for all
$ B \in \mathcal{M}_{\flat }(X), $ then we call that $ P $ is $ \varepsilon -$
admissible. We define $ d_{a}(P) = \inf \{ \varepsilon > 0 : \ P \ \text{is} \
\varepsilon -\text{admissible} \}, $ and call $ d_{a}(P) $ the admissible degree
of $ P. $ We also denote
\begin{align*}
&C_{\varepsilon } (P) = \sup \{ \delta > 0 : \
\widehat{d}_{Y} (P(B), P(B^{\delta })) \leq \varepsilon \
\text{for all} \ B \in \mathcal{M}_{\flat }(X) \}, \\
&C(P) = \inf \{ C_{\varepsilon } (P) : \
P \ \text{is} \ \varepsilon -\text{admissible} \}, \\
&\widehat{h}_{P_{1} \rightarrow P_{2}}(B) =
\sup _{\varepsilon > 0}
\widehat{d}_{Y} (P_{1}(B), P_{2}(B^{\varepsilon })) \quad
\text{for each} \ B \in \mathcal{M}_{\flat }(X), \\
&\widehat{h}(P) = \inf _{B \in \mathcal{M}_{\flat }(X)}
\widehat{h}_{P \rightarrow P}(B), \
\widehat{H}(P) = \sup _{B \in \mathcal{M}_{\flat }(X)}
 \widehat{h}_{P \rightarrow P}(B), \
\text{and} \\
&\widehat{H}(P_{1}, P_{2}) = \sup _{B \in \mathcal{M}_{\flat }(X)}
\max \{\widehat{h}_{P_{1} \rightarrow P_{2}}(B), \
\widehat{h}_{P_{2} \rightarrow P_{1}}(B) \}.
\end{align*}
(2) \ Assume $ P_{1} \neq P_{2}, $ we define
\begin{align*}
O_{\varepsilon }(P_{1}, P_{2}) =
\{\delta > 0 : \ \widehat{d}_{Y} (P_{1}(B), P_{2}(B^{\delta }))
\leq \varepsilon \ \text{and} \ \widehat{d}_{Y} (P_{2}(B), P_{1}(B^{\delta }))
\leq \varepsilon \ \\
\text{for all} \ B \in \mathcal{M}_{\flat } (X) \},
\end{align*}
\begin{align*}
&\omega _{\varepsilon }(P_{1}, P_{2}) = \inf \{\delta : \ \delta \in
O_{\varepsilon }(P_{1}, P_{2}) \}, \quad \Omega _{\varepsilon }(P_{1}, P_{2})
= \sup \{\delta : \ \delta \in
O_{\varepsilon }(P_{1}, P_{2}) \}, \\
&\widehat{\Omega }_{X, Y}(P_{1}, P_{2}) =
\inf \{ \varepsilon > 0 : \ \Omega _{\varepsilon }(P_{1}, P_{2}) = \infty \}.
\end{align*}
We denote $ \omega _{\varepsilon }(P_{1}, P_{2}) =
\Omega _{\varepsilon }(P_{1}, P_{2}) = 0 $
if $ O_{\varepsilon }(P_{1}, P_{2}) = \emptyset , $
and $ \widehat{\Omega }_{X, Y}(P_{1}, P_{2}) = \infty $ if
$ \Omega _{\varepsilon }(P_{1}, P_{2}) = 0 $ for all $ \varepsilon > 0. $
Moreover, $ \widehat{\Omega }_{X, Y}(P_{1}, P_{2}) = 0 $ if
$ P_{1} = P_{2}. $ \\
(3) \ We define
\begin{align*}
O_{\varepsilon }^{\ast }(P_{1} \rightarrow P_{2}) =
\{\delta > 0 : \ \text{for every} \ B \in \mathcal{M}_{\flat } (X),
\ \text{there exist a positive number} \\
\delta _{B} \leq \delta \
\text{such that} \
\widehat{d}_{Y} (P_{1}(B), P_{2}(B^{\delta _{B}})) \leq \varepsilon \},
\end{align*}
$ \widehat{\Omega }_{X, Y}^{\ast }(P_{1}, P_{2}) = \inf \{ \varepsilon > 0 : \
O_{\varepsilon }^{\ast }(P_{1} \rightarrow P_{2}) \ \cap \
O_{\varepsilon }^{\ast }(P_{2} \rightarrow P_{1}) \neq \emptyset \}. $
\par     \vskip  0.2 cm

It is easy to verify that \\
(1) \ $ P $ is $ \varepsilon-$admissible for any $ \varepsilon >
d_{a}(P). $ \\
(2) \ For any $ \lambda \in (0, C(P)) $ and $ \varepsilon \in
(d_{a}(P), \infty ), $ we have \\
$ \widehat{d}_{Y} (P(B), P(B^{\lambda })) < \varepsilon \
(\forall B \in \mathcal{M}_{\flat }(X)). $
\par     \vskip  0.2 cm

{\bf Theorem 4.5.} \ Let $ P, P_{1}, P_{2} \in
\mathfrak{D}_{\flat }(X, Y). $ \\
(1) \ If $ P $ is a Lipschitz map, then $ d_{a}(P) = 0. $ \\
(2) \ $ \max \{ \widehat{\beta }_{X, Y}^{\infty } (P_{1}, P_{2}),
\widehat{H}(P_{1}), \widehat{H}(P_{2}) \} = \widehat{H}(P_{1}, P_{2})
= \max \{ \widehat{\rho } _{s}(P_{1}, P_{2}),
\widehat{H}(P_{1}), \widehat{H}(P_{2}) \}. $ \\
(3) \ If for every $ B \in \mathcal{M}_{\flat }(X), $ the cardinal
$ \sharp \{\widehat{d}_{Y} (P_{1}(B^{\varepsilon }),
P_{2}(B^{\varepsilon })) : \ \varepsilon > 0 \} > 1, $ then \
$ \max \{\widehat{H}(P_{1}), \widehat{H}(P_{2}) \} =
\widehat{H}(P_{1}, P_{2}). $
\par     \vskip  0.1 cm
{\bf Proof.} \ (1) \ Since $ P $ is a Lipschitz map, there exists
a number $ c \geq 0 $ such that
$ \widehat{d}_{Y} (P(B_{1}), P(B_{2})) \leq c \cdot \widehat{d}_{X, H}
(B_{1}, B_{2}) \ (\forall B_{1}, B_{2} \in \mathcal{M}_{\flat }(X)). $
If $ c = 0, $ then obviously $ d_{a}(P) = 0. $
Now we assume $ c > 0. $ For any $ \varepsilon > 0, $ by Lemmas 2.1
and 2.3, it follows easily
that $ \widehat{d}_{X, H}(B, B^{\varepsilon }) \leq \varepsilon , $ so
$ \widehat{d}_{Y} (P(B), P(B^{\varepsilon })) \leq c \cdot
\widehat{d}_{X, H}(B, B^{\varepsilon }) \leq c \varepsilon
\ (\forall B \in \mathcal{M}_{\flat }(X)). $ Hence $ P $
is $ c \varepsilon-$admissible, and so $ d_{a}(P) \leq c \varepsilon , $
which implies $ d_{a}(P) = 0, $ this proves (1). \\
(2) \ We may assume $ P_{1} \neq P_{2}. $ By definition,
one can easily verify that
$ \widehat{\beta }_{X, Y}^{\infty } (P_{1}, P_{2}) \leq
\widehat{H}(P_{1}, P_{2}) $ and $ \widehat{H}(P_{i}) \leq
\widehat{H}(P_{1}, P_{2}) \ (i =1, 2). $ Now let \\
$ \delta > \max \{ \widehat{\beta }_{X, Y}^{\infty } (P_{1}, P_{2}),
\widehat{H}(P_{1}), \widehat{H}(P_{2}) \}, $ then for every $ \lambda
> 0, $ we have
$ \delta > \widehat{\beta }_{X, Y}^{\lambda } (P_{1}, P_{2}), $
so there exists a positive number $ \delta (\lambda ) < \delta $
such that \\
$ \widehat{d}_{Y} (P_{1}(B), P_{2}(B^{\lambda \delta (\lambda )}))
\leq \delta (\lambda ) $ and
$ \widehat{d}_{Y} (P_{2}(B), P_{1}(B^{\lambda \delta (\lambda )}))
\leq \delta (\lambda ) \ (\forall B \in \mathcal{M}_{\flat }(X)). $
Hence for any $ \varepsilon > 0, \\
\widehat{d}_{Y} (P_{1}(B), P_{2}(B^{\varepsilon })) \leq
\max \{ \widehat{d}_{Y} (P_{1}(B),
P_{2}(B^{\lambda \delta (\lambda )})),
\widehat{d}_{Y} (P_{2}(B^{\lambda \delta (\lambda )}),
P_{2}(B^{\varepsilon })) \} \\
\leq \max \{ \delta , \widehat{H}(P_{2}) \}, $
so $ \widehat{h}_{P_{1} \rightarrow P_{2}}(B)
\leq \max \{ \delta , \widehat{H}(P_{2}) \}. $ \\
Similarly, $ \widehat{h}_{P_{2} \rightarrow P_{1}}(B)
\leq \max \{ \delta , \widehat{H}(P_{1}) \}. $ \\
Therefore $ \widehat{H}(P_{1}, P_{2}) \leq
\max \{ \delta , \widehat{H}(P_{1}), \widehat{H}(P_{2}) \}
= \delta , $ which implies \\
$ \widehat{H}(P_{1}, P_{2}) \leq
\max \{ \widehat{\beta }_{X, Y}^{\infty } (P_{1}, P_{2}),
\widehat{H}(P_{1}), \widehat{H}(P_{2}) \}, $ and the first equality
holds. \\
For the second equality, by the above discussion, we only need to
show that $ \widehat{H}(P_{1}, P_{2})
\leq \max \{ \widehat{\rho } _{s}(P_{1}, P_{2}),
\widehat{H}(P_{1}), \widehat{H}(P_{2}) \}. $ To see this, let
$ \delta > \max \{ \widehat{\rho } _{s}(P_{1}, P_{2}),
\widehat{H}(P_{1}), \widehat{H}(P_{2}) \}, $ then for every
$ B \in \mathcal{M}_{\flat }(X)), $ we have \\
$ \delta > \widehat{d}_{Y} (P_{1}(B), P_{2}(B)) $ \ and \
$ \delta > \widehat{d}_{Y} (P_{i}(B),
P_{i}(B^{\varepsilon })) $ for all $ \varepsilon > 0
\ (i = 1, 2). $ \\
So $ \widehat{d}_{Y} (P_{1}(B), P_{2}(B^{\varepsilon }))
\leq \max \{ \widehat{d}_{Y} (P_{1}(B), P_{2}(B)),
\widehat{d}_{Y} (P_{2}(B), P_{2}(B^{\varepsilon })) \} <
\delta . $ \\
Hence $ \widehat{h}_{P_{1} \rightarrow P_{2}}(B) \leq \delta . $
Similarly, $ \widehat{h}_{P_{2} \rightarrow P_{1}}(B) \leq \delta . $
Therefore $ \widehat{H}(P_{1}, P_{2}) \leq \delta , $ which
implies $ \widehat{H}(P_{1}, P_{2}) \leq
\max \{ \widehat{\rho } _{s}(P_{1}, P_{2}),
\widehat{H}(P_{1}), \widehat{H}(P_{2}) \}. $ This proves (2). \\
(3) \ For every $ B \in \mathcal{M}_{\flat }(X), $ by assumption,
there exists an $ \varepsilon > 0 $ such that
$ \widehat{d}_{Y} (P_{1}(B^{\varepsilon }), P_{2}(B^{\varepsilon }))
\neq \widehat{d}_{Y} (P_{1}(B), P_{2}(B)). $ We may as well assume
that \\
$ \widehat{d}_{Y} (P_{1}(B^{\varepsilon }), P_{2}(B^{\varepsilon }))
< \widehat{d}_{Y} (P_{1}(B), P_{2}(B)). $ Then by the inequality \\
$ \widehat{d}_{Y} (P_{1}(B), P_{2}(B)) \leq
\max \{ \widehat{d}_{Y} (P_{1}(B), P_{1}(B^{\varepsilon })),
\widehat{d}_{Y} (P_{1}(B^{\varepsilon }), P_{2}(B^{\varepsilon })),
\widehat{d}_{Y} (P_{2}(B^{\varepsilon }), P_{2}(B)) \}, $ we get \\
$ \widehat{d}_{Y} (P_{1}(B), P_{2}(B)) \leq
\max \{ \widehat{d}_{Y} (P_{1}(B), P_{1}(B^{\varepsilon })),
\widehat{d}_{Y} (P_{2}(B), P_{2}(B^{\varepsilon })) \} \\
\leq \max \{\widehat{h}_{P_{1} \rightarrow P_{1}}(B),
\widehat{h}_{P_{2} \rightarrow P_{2}}(B) \}. $
Therefore \\
$ \widehat{\rho } _{s}(P_{1}, P_{2}) \leq \sup _{B \in
\mathcal{M}_{\flat }(X))} \max \{\widehat{h}_{P_{1} \rightarrow P_{1}}(B),
\widehat{h}_{P_{2} \rightarrow P_{2}}(B) \} =
\max \{\widehat{H}(P_{1}), \widehat{H}(P_{2}) \}. $
Now we may as well assume that
$ \widehat{H}(P_{1}) \leq \widehat{H}(P_{2}). $ Then for
any $ \delta > 0, $ we have \\
$ \widehat{d}_{Y} (P_{1}(B), P_{2}(B^{\delta })) \leq
\max \{ \widehat{d}_{Y} (P_{1}(B), P_{2}(B)),
\widehat{d}_{Y} (P_{2}(B), P_{2}(B^{\delta })) \} \\
\leq \max \{ \widehat{\rho } _{s}(P_{1}, P_{2}),
\widehat{H}(P_{2}) \} = \widehat{H}(P_{2}), $ \
similarly, $ \widehat{d}_{Y} (P_{2}(B), P_{1}(B^{\delta }))
\leq \widehat{H}(P_{2}), $ which implies
$ \widehat{H}(P_{1}, P_{2}) \leq \widehat{H}(P_{2}), $
so the equality holds. This proves (3), and the proof of
Theorem 4.5 is completed. \quad $ \Box $
\par     \vskip  0.2 cm

{\bf Theorem 4.6.} \ Let $ P_{1}, P_{2} \in
\mathfrak{D}_{\flat }(X, Y) $ with $ P_{1} \neq P_{2}. $ We have \\
(1) \ $ \widehat{\beta }_{X, Y}^{0} (P_{1}, P_{2}) =
\widehat{\rho } _{s}(P_{1}, P_{2}) \ \Longleftrightarrow
\ \max \{ d_{a}(P_{1}), d_{a}(P_{2}) \} \leq
\widehat{\rho } _{s}(P_{1}, P_{2}). $ \\
(2) \ $ \widehat{\beta }_{X, Y}^{\infty } (P_{1}, P_{2}) =
\widehat{\Omega }_{X, Y}(P_{1}, P_{2}). $ \\
(3) \ $ \widehat{\beta }_{X, Y}^{\ast 0}(P_{1}, P_{2}) =
\widehat{\rho } _{s}(P_{1}, P_{2}). $ \\
(4) \ $ \widehat{\beta }_{X, Y}^{\ast \infty }(P_{1}, P_{2})
= \widehat{\Omega }_{X, Y}^{\ast }(P_{1}, P_{2}). $
\par     \vskip  0.1 cm
{\bf Proof.} \ (1) \ If $ \widehat{\beta }_{X, Y}^{0} (P_{1}, P_{2}) =
\widehat{\rho } _{s}(P_{1}, P_{2}), $ then for any $ \varepsilon > 0, $
there exists an $ \delta > 0 $ such that for any $ \lambda \in (0, \delta ), $
we have $ \widehat{\beta }_{X, Y}^{\lambda } (P_{1}, P_{2}) <
\widehat{\rho } _{s}(P_{1}, P_{2}) + \varepsilon . $ By definition,
there exists a positive number $ \varepsilon (\lambda ) : \
\widehat{\beta }_{X, Y}^{\lambda } (P_{1}, P_{2}) \leq
\varepsilon (\lambda ) < \widehat{\rho } _{s}(P_{1}, P_{2}) + \varepsilon
$ such that
$ \widehat{d}_{Y} (P_{1}(B), P_{2}(B^{\lambda \varepsilon (\lambda ) }))
\leq \varepsilon (\lambda ) $ and
$ \widehat{d}_{Y} (P_{2}(B), P_{1}(B^{\lambda \varepsilon (\lambda ) }))
\leq \varepsilon (\lambda ) \ ( \forall B \in \mathcal{M}_{\flat } (X)). $
So $ \widehat{d}_{Y} (P_{1}(B), P_{1}(B^{\lambda \varepsilon (\lambda ) }))
\leq \max \{ \widehat{d}_{Y} (P_{1}(B), P_{2}(B^{\lambda \varepsilon (\lambda ) })),
\widehat{d}_{Y} (P_{2}(B^{\lambda \varepsilon (\lambda ) }),
P_{1}(B^{\lambda \varepsilon (\lambda ) })) \} \leq
\varepsilon (\lambda ). $
Similarly, $ \widehat{d}_{Y} (P_{2}(B), P_{2}(B^{\lambda \varepsilon (\lambda ) }))
\leq \varepsilon (\lambda ). $ So both $ P_{1} $ and $ P_{2} $ are
$ \varepsilon (\lambda )-$admissible. Hence by definition, we get
$ \max \{ d_{a}(P_{1}), d_{a}(P_{2}) \} \leq \varepsilon (\lambda )
< \widehat{\rho } _{s}(P_{1}, P_{2}) + \varepsilon , $ which implies
$ \max \{ d_{a}(P_{1}), d_{a}(P_{2}) \} \leq
\widehat{\rho } _{s}(P_{1}, P_{2}) $ by the arbitrary choice of
$ \varepsilon . $ \\
Conversely, if $ \max \{ d_{a}(P_{1}), d_{a}(P_{2}) \} \leq
\widehat{\rho } _{s}(P_{1}, P_{2}), $ then we need to show that \\
$ \widehat{\beta }_{X, Y}^{0} (P_{1}, P_{2}) =
\widehat{\rho } _{s}(P_{1}, P_{2}). $ If otherwise, then
$ \widehat{\beta }_{X, Y}^{0} (P_{1}, P_{2}) >
\widehat{\rho } _{s}(P_{1}, P_{2}). $ Take $ \varepsilon _{0} =
\frac{1}{2} (\widehat{\beta }_{X, Y}^{0} (P_{1}, P_{2}) +
\widehat{\rho } _{s}(P_{1}, P_{2}) ), $ then
$ \widehat{\rho } _{s}(P_{1}, P_{2}) < \varepsilon _{0} <
\widehat{\beta }_{X, Y}^{0} (P_{1}, P_{2}). $ So for every $ \lambda > 0, \
 \widehat{\beta }_{X, Y}^{\lambda } (P_{1}, P_{2}) > \varepsilon _{0}, $
and then there exists a $ B_{\lambda } \in \mathcal{M}_{\flat } (X) $
such that \\
$ \widehat{d}_{Y} (P_{1}(B_{\lambda }),
P_{2}(B_{\lambda }^{\lambda \varepsilon _{0}})) > \varepsilon _{0} $ or
$ \widehat{d}_{Y} (P_{2}(B_{\lambda }),
P_{1}(B_{\lambda }^{\lambda \varepsilon _{0}})) > \varepsilon _{0}. $ \\
On the other hand, we have $ d_{a}(P_{i }) < \varepsilon _{0},  \ i =1, 2. $
So by definition, there exist positive numbers $ \varepsilon _{i} <
\varepsilon _{0} $ such that $ P_{i} $ are
$ \varepsilon _{i}-$admissible, i.e., there exist positive numbers
$ \delta _{i} $ such that \\
$ \widehat{d}_{Y} (P_{i}(B), P_{i}(B^{\delta _{i}})) \leq
\varepsilon _{i} \ (i =1, 2) \
(\forall B \in \mathcal{M}_{\flat } (X)). $ \\
Now take a $ \lambda _{0} > 0 $ such that
$ \lambda _{0} \varepsilon _{0} < \min \{ \delta _{1},
\delta _{2} \}. $ Then by Lemma 2.3,
$ (B^{\lambda _{0} \varepsilon _{0}})^{\delta _{i}} = B^{\delta _{i}}, $
so \\
$ \widehat{d}_{Y} (P_{i}(B), P_{i}(B^{\lambda _{0} \varepsilon _{0}}))
\leq \max \{\widehat{d}_{Y} (P_{i}(B), P_{i}(B^{\delta _{i}})),
\widehat{d}_{Y} (P_{i}((B^{\lambda _{0} \varepsilon _{0}})^{\delta _{i}}),
P_{i}(B^{\lambda _{0} \varepsilon _{0}}))
\} \leq \varepsilon _{i} \\
(i =1, 2) \ (\forall B \in \mathcal{M}_{\flat } (X)). $ \\
In particular,
$ \widehat{d}_{Y} (P_{i}(B_{\lambda _{0}}),
P_{i}(B_{\lambda _{0}}^{\lambda _{0} \varepsilon _{0}})) \leq \varepsilon _{i}
< \varepsilon _{0} \ (i =1, 2). $ \\
Since
$ \widehat{d}_{Y} (P_{1}(B_{\lambda _{0}}^{\lambda _{0} \varepsilon _{0}}),
P_{2}(B_{\lambda _{0}}^{\lambda _{0} \varepsilon _{0}})) \leq
\widehat{\rho } _{s}(P_{1}, P_{2}) < \varepsilon _{0}, $ we get \\
$ \widehat{d}_{Y} (P_{1}(B_{\lambda _{0}}),
P_{2}(B_{\lambda _{0}}^{\lambda _{0} \varepsilon _{0}}))
\leq \max \{\widehat{d}_{Y} (P_{1}(B_{\lambda _{0}}),
P_{1}(B_{\lambda _{0}}^{\lambda _{0} \varepsilon _{0}})),
\widehat{d}_{Y} (P_{1}(B_{\lambda _{0}}^{\lambda _{0} \varepsilon _{0}}),
P_{2}(B_{\lambda _{0}}^{\lambda _{0} \varepsilon _{0}}))
\} < \varepsilon _{0}, $ similarly,
$ \widehat{d}_{Y} (P_{2}(B_{\lambda _{0}}),
P_{1}(B_{\lambda _{0}}^{\lambda _{0} \varepsilon _{0}}))
< \varepsilon _{0}, $ a contradiction! Hence we must have
$ \widehat{\beta }_{X, Y}^{0} (P_{1}, P_{2}) =
\widehat{\rho } _{s}(P_{1}, P_{2}). $ This proves (1). \\
(2) \ For any $ \varepsilon > \widehat{\beta }_{X, Y}^{\infty } (P_{1}, P_{2}), $
by Lemma 4.3, $ \varepsilon > \widehat{\beta }_{X, Y}^{\lambda } (P_{1}, P_{2})
\ (\forall \lambda > 0). $ So for every $ \lambda > 0, $ by definition,
there exists a positive number $ \varepsilon (\lambda ) < \varepsilon $
such that \\
$ \widehat{d}_{Y} (P_{1}(B), P_{2}(B^{\lambda \varepsilon (\lambda )}))
\leq \varepsilon (\lambda ) $ and
$ \widehat{d}_{Y} (P_{2}(B), P_{1}(B^{\lambda \varepsilon (\lambda )}))
\leq \varepsilon (\lambda ) \
(\forall B \in \mathcal{M}_{\flat } (X)). $ \\
So $ \lambda \varepsilon (\lambda )
\in O_{\varepsilon (\lambda ) }(P_{1}, P_{2}) \subset
O_{\varepsilon }(P_{1}, P_{2}). $ Then by Theorem 4.2.(2), we have
$ \varepsilon (\lambda ) \geq \widehat{\beta }_{X, Y}^{\lambda } (P_{1}, P_{2})
\geq \widehat{\rho } _{s}(P_{1}, P_{2}). $ Since $ P_{1} \neq P_{2}, $
it is easy to see that $ \widehat{\rho } _{s}(P_{1}, P_{2}) > 0. $
Hence $ \lambda \varepsilon (\lambda ) \rightarrow \infty $ as
$ \lambda \rightarrow \infty , $ which implies that
$ \Omega _{\varepsilon }(P_{1}, P_{2}) = \infty . $ Then by definition,
$ \widehat{\Omega }_{X, Y}(P_{1}, P_{2}) \leq \varepsilon , $ and
so $ \widehat{\Omega }_{X, Y}(P_{1}, P_{2}) \leq
\widehat{\beta }_{X, Y}^{\infty } (P_{1}, P_{2}) $ by the arbitrary
choice of $ \varepsilon . $ \\
Conversely, for any $ \varepsilon > \widehat{\Omega }_{X, Y}(P_{1}, P_{2}), $
by definition, there exists a positive number $ \varepsilon ^{\prime }
< \varepsilon $ such that
$ \Omega _{\varepsilon ^{\prime } }(P_{1}, P_{2}) = \infty , $
in particular, $ O_{\varepsilon ^{\prime }}(P_{1}, P_{2}) \neq \emptyset . $
Then, for any $ \lambda > 0, $ there exists an $ \delta \in
O_{\varepsilon ^{\prime }}(P_{1}, P_{2}) \cap (\lambda \varepsilon , \infty ). $
By definition, \\
$ \widehat{d}_{Y} (P_{1}(B), P_{2}(B^{\delta })) \leq \varepsilon ^{\prime } $
and $ \widehat{d}_{Y} (P_{2}(B), P_{1}(B^{\delta })) \leq
\varepsilon ^{\prime } \ (\forall B \in \mathcal{M}_{\flat } (X)). $ \\
Let $ \lambda ^{\prime } = \delta / \varepsilon , $ then $ \lambda ^{\prime }
> \lambda , $ and \\
$ \widehat{d}_{Y} (P_{1}(B), P_{2}(B^{\lambda ^{\prime } \varepsilon })) <
\varepsilon $ and $ \widehat{d}_{Y} (P_{2}(B),
P_{1}(B^{\lambda ^{\prime } \varepsilon }))
< \varepsilon \ (\forall B \in \mathcal{M}_{\flat } (X)). $ \\
So $ \widehat{\beta }_{X, Y}^{\lambda ^{\prime } } (P_{1}, P_{2})
\leq \varepsilon , $ and then by Lemma 4.3,
$ \widehat{\beta }_{X, Y}^{\lambda } (P_{1}, P_{2})
\leq \varepsilon , $ which implies
$ \widehat{\beta }_{X, Y}^{\infty } (P_{1}, P_{2})
\leq \varepsilon . $ Therefore,
$ \widehat{\beta }_{X, Y}^{\infty } (P_{1}, P_{2}) \leq
\widehat{\Omega }_{X, Y}(P_{1}, P_{2}), $ and the equality holds,
this proves (2). \\
(3) \ By Theorem 4.2.(2), $ \widehat{\beta }_{X, Y}^{\ast 0}(P_{1}, P_{2})
\leq \widehat{\rho } _{s}(P_{1}, P_{2}). $
As in the proof of (2) above, $ \widehat{\rho } _{s}(P_{1}, P_{2}) > 0. $
For any positive number $ \varepsilon < \widehat{\rho } _{s}(P_{1}, P_{2}), $
by definition, there exists a $ B_{0} \in \mathcal{M}_{\flat } (X) $
such that $ \widehat{d}_{Y} (P_{1}(B_{0}), P_{2}(B_{0})) > \varepsilon . $
Denote $ \delta = \text{diam} B_{0} / (2 \widehat{\rho } _{s}(P_{1}, P_{2}))
> 0, $ and take a $ \lambda \in (0, \delta ). $ Then for any $ \varepsilon _{1} \in
(\widehat{\beta }_{X, Y}^{\ast \lambda}(P_{1}, P_{2}),
2 \widehat{\rho } _{s}(P_{1}, P_{2})), $ there exists an
$ \varepsilon _{1}^{\prime } \in (0, \varepsilon _{1}) $ such that \\
$ \widehat{d}_{Y} (P_{1}(B), P_{2}(B^{\lambda \varepsilon _{B}(\lambda )}))
\leq \varepsilon _{1}^{\prime } $ \ and \
$ \widehat{d}_{Y}
(P_{2}(B), P_{1}(B^{\lambda \varepsilon _{B}^{\prime }(\lambda )}))
\leq \varepsilon _{1}^{\prime } $ \\
with some $ \varepsilon _{B}(\lambda ),
\varepsilon _{B}^{\prime }(\lambda ) \in (0, \varepsilon _{1}^{\prime } ] $
for every $ B \in \mathcal{M}_{\flat } (X). $ In particular, \\
$ \widehat{d}_{Y} (P_{1}(B_{0}),
P_{2}(B_{0}^{\lambda \varepsilon _{B_{0}}(\lambda )}))
\leq \varepsilon _{1}^{\prime } $ \ and \
$ \widehat{d}_{Y}
(P_{2}(B_{0}), P_{1}(B_{0}^{\lambda \varepsilon _{B_{0}}^{\prime }(\lambda )}))
\leq \varepsilon _{1}^{\prime }. $ \\
Since $ \lambda \varepsilon _{B_{0}}(\lambda ) \leq \lambda
\varepsilon _{1}^{\prime } < \lambda \varepsilon _{1} <
2 \lambda \widehat{\rho } _{s}(P_{1}, P_{2}) < 2 \delta
\widehat{\rho } _{s}(P_{1}, P_{2}) = \text{diam} B_{0}, $ by Lemma 2.3,
$ B_{0}^{\lambda \varepsilon _{B_{0}}(\lambda )} = B_{0}, $ so
$ \widehat{d}_{Y} (P_{1}(B_{0}), P_{2}(B_{0})) =
\widehat{d}_{Y} (P_{1}(B_{0}),
P_{2}(B_{0}^{\lambda \varepsilon _{B_{0}}(\lambda )}))
\leq \varepsilon _{1}^{\prime } < \varepsilon _{1}. $ Since
$ \varepsilon _{1} $ is arbitrary chosen, we get
$ \widehat{d}_{Y} (P_{1}(B_{0}), P_{2}(B_{0})) \leq
\widehat{\beta }_{X, Y}^{\ast \lambda}(P_{1}, P_{2}), $
hence $ \varepsilon <
\widehat{\beta }_{X, Y}^{\ast \lambda}(P_{1}, P_{2})
\leq \widehat{\beta }_{X, Y}^{\ast 0}(P_{1}, P_{2}), $ which implies
$ \widehat{\rho } _{s}(P_{1}, P_{2}) \leq
\widehat{\beta }_{X, Y}^{\ast 0}(P_{1}, P_{2}), $ and the equality
holds, This proves (3). \\
(4) \ For any $ \varepsilon >
\widehat{\beta }_{X, Y}^{\ast \infty }(P_{1}, P_{2}), $ there exists
an $ \lambda > 0 $ such that
$ \widehat{\beta }_{X, Y}^{\ast \lambda }(P_{1}, P_{2}) < \varepsilon . $
By definition, there exists an $ \varepsilon ^{\prime } \in
(0, \varepsilon ) $ such that \\
$ \widehat{d}_{Y} (P_{1}(B), P_{2}(B^{\lambda \varepsilon _{B}}))
\leq \varepsilon ^{\prime } $ and
$ \widehat{d}_{Y} (P_{2}(B), P_{1}(B^{\lambda \varepsilon _{B}^{\prime }}))
\leq \varepsilon ^{\prime } $ \\
with some $ \varepsilon _{B},
\varepsilon _{B}^{\prime } \in (0, \varepsilon ^{\prime }] $ for
every $ B \in \mathcal{M}_{\flat } (X). $ So
$ \lambda \varepsilon ^{\prime } \in
O_{\varepsilon ^{\prime } }^{\ast }(P_{1} \rightarrow P_{2}) \ \cap \
O_{\varepsilon ^{\prime } }^{\ast }(P_{2} \rightarrow P_{1}), $ and
then $ \widehat{\Omega }_{X, Y}^{\ast }(P_{1}, P_{2}) \leq
\varepsilon ^{\prime } < \varepsilon , $ which implies
$ \widehat{\Omega }_{X, Y}^{\ast }(P_{1}, P_{2}) \leq
\widehat{\beta }_{X, Y}^{\ast \infty }(P_{1}, P_{2}). $ \\
Conversely, for any $ \varepsilon >
\widehat{\Omega }_{X, Y}^{\ast }(P_{1}, P_{2}), $ there exists an
$ \varepsilon ^{\prime } \in (0, \varepsilon ) $ such that
$ O_{\varepsilon ^{\prime } }^{\ast }(P_{1} \rightarrow P_{2}) \ \cap \
O_{\varepsilon ^{\prime } }^{\ast }(P_{2} \rightarrow P_{1}) \neq
\emptyset . $ Let $ \delta $ be a positive number in this intersection
set, then, for every $ B \in \mathcal{M}_{\flat } (X), $ there exist
$ \delta _{B}, \delta _{B}^{\prime } \in (0, \delta ] $ such that \\
$ \widehat{d}_{Y} (P_{1}(B), P_{2}(B^{\delta _{B}}))
\leq \varepsilon ^{\prime } $ and
$ \widehat{d}_{Y} (P_{2}(B), P_{1}(B^{\delta _{B}^{\prime }}))
\leq \varepsilon ^{\prime }. $ \\
Take $ \lambda _{0} = \delta / \varepsilon ^{\prime }, \
\varepsilon _{B} = \delta _{B} / \lambda _{0} $ and
$ \varepsilon _{B}^{\prime } = \delta _{B}^{\prime } / \lambda _{0}, $
then $ \lambda _{0} > 0, \ \varepsilon _{B} \leq \varepsilon ^{\prime },
\ \varepsilon _{B}^{\prime } \leq \varepsilon ^{\prime }, $ and from
$ \widehat{d}_{Y} (P_{1}(B), P_{2}(B^{\lambda _{0} \varepsilon _{B}}))
\leq \varepsilon ^{\prime } $ and $ \widehat{d}_{Y} (P_{2}(B),
P_{1}(B^{\lambda _{0} \varepsilon _{B}^{\prime }}))
\leq \varepsilon ^{\prime } \
(\forall B \in \mathcal{M}_{\flat } (X)), $ we get
$ \widehat{\beta }_{X, Y}^{\ast \infty }(P_{1}, P_{2}) \leq
\widehat{\beta }_{X, Y}^{\ast \lambda _{0}}(P_{1}, P_{2}) \leq
\varepsilon ^{\prime } < \varepsilon , $ which implies
$ \widehat{\beta }_{X, Y}^{\ast \infty }(P_{1}, P_{2}) \leq
\widehat{\Omega }_{X, Y}^{\ast }(P_{1}, P_{2}), $ and the
equality holds. This proves (4), and the proof of Theorem 4.6
is completed.  \quad $ \Box $
\par     \vskip  0.2 cm

{\bf Corollary 4.7.} \ Let $ P_{1}, P_{2} \in
\mathfrak{D}_{\flat }(X, Y). $ \\
(1) \ If both $ P_{1} $ and $ P_{2} $ are Lipschitz maps, then
$ \widehat{\beta }_{X, Y}^{0} (P_{1}, P_{2}) =
\widehat{\rho } _{s}(P_{1}, P_{2}). $ \\
(2) \ If $ \max \{ \widehat{H}(P_{1}), \widehat{H}(P_{2}) \}
< \widehat{H}(P_{1}, P_{2}), $ then for every $ \lambda > 0, $
we have \\
$ \widehat{\beta }_{X, Y}^{\lambda }(P_{1}, P_{2}) =
\widehat{\beta }_{X, Y}^{0} (P_{1}, P_{2}) =
\widehat{\rho } _{s}(P_{1}, P_{2}) =
\widehat{H}(P_{1}, P_{2}). $ \\
(3) \ If $ \max \{ d_{a}(P_{1}), d_{a}(P_{2}) \}
> \widehat{\rho } _{s}(P_{1}, P_{2}), $ then
$ \max \{ \widehat{H}(P_{1}), \widehat{H}(P_{2}) \} =
\widehat{H}(P_{1}, P_{2}). $
\par     \vskip  0.1 cm
{\bf Proof.} \ (1) follows from Theorem 4.5.(1)
and Theorem 4.6.(1). \\
(2) follows from Lemma 4.3 and Theorem 4.5.(2). \\
(3) follows from Theorem 4.5.(2) and
Theorem 4.6.(1).  \quad $ \Box $
\par     \vskip  0.2 cm

{\bf Example 4.8.} \ Let $ X = \C_{p} $ be the Tate field with
the canonical $ p-$adic metric $ \widehat{d}_{p} $ as above.
It is well known that $ \mid \C_{p} \mid _{p} = p^{\Q} $
(see [Sc], p.45), and it follows
that $ \text{diam}B_{a}(r) = r $ for any $ a \in \C_{p} $ and
$ r > 0. $ Now we define two maps
$ P_{i} : \ \mathcal{M}_{\flat } (X) \longrightarrow
\mathcal{M}_{\flat } (X) \ (i = 1, 2) $ as follows: \\
We define $ P_{1}(B) = B^{1 / \text{diam}B} $ for all $ B \in
\mathcal{M}_{\flat } (X); $ \ and \
$ P_{2}(B) = P_{1}(B) $ if $ B \neq B_{0}(1), $
while $ P_{2}(B_{0}(1)) = B_{0}(2). $ \\
We assert that $ d_{a}(P_{1}) = \infty , $ that is, for every
$ \varepsilon > 0, P_{1} $ is not $ \varepsilon-$admissible.
To see this, for any $ \delta > 0, $ take a positive number
$ r < \min \{\delta , \frac{1}{\delta },
\frac{1}{\varepsilon } \}, $ denote
$ \delta _{0} = \max \{\delta , \frac{1}{\delta } \} $
and let $ B = B_{0}(r). $ Then
by Lemma 2.3, $ P_{1}(B) = P_{1}(B_{0}(r)) = B^{1 / r} =
B_{0}(1 / r) $ and $ P_{1}(B^{\delta }) = P_{1}(B_{0}(\delta ))
= B_{0}(\delta )^{1 / \delta } = B_{0}( \delta _{0}). $
So by Lemma 2.1, we have \\
$ \widehat{d}_{p, H} (P_{1}(B), P_{1}(B^{\delta })) =
\widehat{d}_{p, H} (B_{0}(1 / r), B_{0}( \delta _{0})) = 1 / r
> \varepsilon . $ Therefore, $ d_{a}(P_{1}) = \infty . $ On the
other hand, by the above definition and Lemma 2.1, we have
$ \widehat{\rho } _{s}(P_{1}, P_{2}) =
\sup _{B \in \mathcal{M}_{\flat } (X)}
\widehat{d}_{p, H} (P_{1}(B), P_{2}(B)) =
\widehat{d}_{p, H} (P_{1}(B_{0}(1)), P_{2}(B_{0}(1))) =
\widehat{d}_{p, H}(B_{0}(1), B_{0}(2)) = 2. $ Hence
$ d_{a}(P_{1}) > \widehat{\rho } _{s}(P_{1}, P_{2}), $ so by
Theorem 4.6.(1), we get
$ \widehat{\beta }_{X, Y}^{0}(P_{1}, P_{2}) >
\widehat{\rho } _{s}(P_{1}, P_{2}) $ with
$ Y = \mathcal{M}_{\flat } (X) = \mathcal{M}_{\flat } (\C_{p}).
\quad  \Box $
\par     \vskip  0.2 cm

Now let $ K $ be a complete non-Archimedean valued field with
absolute value $ | \ |, $ and $ X $ be a non-Archimedean normed
linear space over $ K $ with norm $ \| \ \| $ (see [FP], chapter 1).
As before, $ K^{\times } = K \backslash \{0 \}. $ Let
$ (Y, \widehat{d}_{Y} ) $ be a non-Archimedean metric space. For
$ a \in K^{\times } $ and $ B \in \mathcal{M}_{\flat } (X), $
it is easy to see that $ a \cdot B \in \mathcal{M}_{\flat } (X), $
more precisely, $ a \cdot B_{b}(r) = B_{a b}(| a | r) $ and
$ a \cdot \overline{B}_{b}(r) = \overline{B}_{a b}(| a | r) \
(b \in X \ \text{and} \ r > 0). $ Moreover, for any
$ \varepsilon > 0, $ we have
$ (a \cdot B)^{\varepsilon } = a \cdot B^{\varepsilon / | a |}. $
Here $ a \cdot B = \{a \cdot v : \ v \in B \}. $
For $ P \in \mathfrak{D}_{\flat }(X, Y), $ we define
$ P^{a}(B) = P(a \cdot B), $ and when $ Y = K, $
we also define $ (a P)(B) = a \cdot P(B). $ Obviously,
$ P^{a} \in \mathfrak{D}_{\flat }(X, Y). $ When
$ Y = K, $ we also have
$ a P \in \mathfrak{D}_{\flat }(X, Y). $
\par     \vskip  0.2 cm

{\bf Theorem 4.9.} \ let $ (K, | \ |) $ be a complete
non-Archimedean valued field, $ (X, \| \ \|) $ be a
non-Archimedean normed linear space over $ K, $ and
$ (Y, \widehat{d}_{Y} ) $ be a non-Archimedean metric space.
For any $ P_{1}, P_{2} \in \mathfrak{D}_{\flat }(X, Y) $
and $ a \in K^{\times }, $ we have \\
$ \widehat{\beta }_{X, Y}^{| a |}(P_{1}, P_{2}) =
\widehat{\beta }_{X, Y}(P_{1}^{a}, P_{2}^{a}) $ \ and \
$ \widehat{\beta }_{X, Y}^{\ast | a |}(P_{1}, P_{2}) =
\widehat{\beta }_{X, Y}^{\ast }(P_{1}^{a}, P_{2}^{a}). $ \\
In particular, if $ Y = K, $ then we have \\
$ \widehat{\beta }_{X, Y}^{| a |}(P_{1}, P_{2}) =
| a |^{-1} \cdot \widehat{\beta }_{X, Y}(a P_{1}, a P_{2}) $ \
and \ $ \widehat{\beta }_{X, Y}^{\ast | a |}(P_{1}, P_{2}) =
| a |^{-1} \cdot
\widehat{\beta }_{X, Y}^{\ast }(a P_{1}, a P_{2}). $
\par     \vskip  0.1 cm
{\bf Proof.} \ For the first equality, let $ \varepsilon >
\widehat{\beta }_{X, Y}^{| a |}(P_{1}, P_{2}). $ Then
there exists a positive number $ \varepsilon ^{\prime } <
\varepsilon $ such that \\
$ \widehat{d}_{Y} (P_{1}(B), P_{2}(B^{| a | \varepsilon ^{\prime }}))
\leq \varepsilon ^{\prime } $ and
$ \widehat{d}_{Y} (P_{2}(B), P_{1}(B^{| a | \varepsilon ^{\prime }}))
\leq \varepsilon ^{\prime } \
(\forall B \in \mathcal{M}_{\flat } (X)). $ So \\
$ \widehat{d}_{Y} (P_{1}^{a}(B), P_{2}^{a}(B^{ \varepsilon ^{\prime }}))
= \widehat{d}_{Y} (P_{1}(a \cdot B),
P_{2}((a \cdot B)^{| a | \varepsilon ^{\prime }}))
\leq \varepsilon ^{\prime }, $ \\
similarly,
$ \widehat{d}_{Y} (P_{2}^{a}(B), P_{1}^{a}(B^{ \varepsilon ^{\prime }}))
\leq \varepsilon ^{\prime }, $ which implies
$ \widehat{\beta }_{X, Y}(P_{1}^{a}, P_{2}^{a}) \leq
\varepsilon ^{\prime } < \varepsilon , $ and so
$ \widehat{\beta }_{X, Y}(P_{1}^{a}, P_{2}^{a}) \leq
\widehat{\beta }_{X, Y}^{| a |}(P_{1}, P_{2}). $ \\
Conversely, let $ \varepsilon >
\widehat{\beta }_{X, Y}(P_{1}^{a}, P_{2}^{a}), $ then there
exists a positive number $ \varepsilon ^{\prime } <
\varepsilon $ such that \\
$ \widehat{d}_{Y} (P_{1}^{a}(B), P_{2}^{a}(B^{ \varepsilon ^{\prime }}))
\leq \varepsilon ^{\prime } $ and
$ \widehat{d}_{Y} (P_{2}^{a}(B), P_{1}^{a}(B^{ \varepsilon ^{\prime }}))
\leq \varepsilon ^{\prime } \ (\forall B \in \mathcal{M}_{\flat } (X)). $
So \\
$ \widehat{d}_{Y} (P_{1}(B), P_{2}(B^{| a | \varepsilon ^{\prime }}))
= \widehat{d}_{Y} (P_{1}^{a}(a^{-1} B),
P_{2}^{a}((a^{-1} B)^{ \varepsilon ^{\prime }}))
\leq \varepsilon ^{\prime }, $ \\
similarly,
$ \widehat{d}_{Y} (P_{2}(B), P_{1}(B^{| a | \varepsilon ^{\prime }}))
\leq \varepsilon ^{\prime }, $ which implies
$ \widehat{\beta }_{X, Y}^{| a |}(P_{1}, P_{2}) \leq
\varepsilon ^{\prime } < \varepsilon , $ and so
$ \widehat{\beta }_{X, Y}^{| a |}(P_{1}, P_{2})
\leq \widehat{\beta }_{X, Y}(P_{1}^{a}, P_{2}^{a}). $ So the first
equality holds. The second equality can be similarly done.  \\
For the case $ Y = K, $ note that
$ \widehat{d}_{Y} ((a P_{1})(B), (a P_{2})(B)) = | a | \cdot
| P_{1}(B)- P_{2}(B)|, $ so \\
$ \widehat{d}_{Y} ((a P_{1})(B), (a P_{2})(B^{ \varepsilon }))
\leq \varepsilon \Longleftrightarrow
\widehat{d}_{Y} (P_{1}(B),
P_{2}(B^{| a | \cdot \varepsilon / | a |})) \leq
\varepsilon / | a |. $ Then the other two equalities follows
easily by the definition. The proof is completed.
\quad $ \Box $
\par     \vskip  0.2 cm

Recall that a non-Archimedean valued ring is a commutative ring
$ A $ with a non-Archimedean absolute value $ | \ |, $
i.e., a function $ | \ |: \ A \rightarrow \R $ satisfying
the rules: \\
1. $ | a | \geq 0 $ and $ | a | = 0 \Leftrightarrow a = 0. $ \
2. $ | a b | = | a | \cdot | b |. $ \ 3. $ | a + b | \leq
\max \{| a |, \ | b | \}. $ \\
If $ A = K $ is a field, then $ (K, | \ |) $ is a
non-Archimedean valued field as before (see [BGR] and [FP]). Now let
$ (X, \widehat{d}_{X} ) $ be a non-Archimedean metric space. For
$ f, g \in \mathfrak{M}(X \rightarrow A) $ and $ a \in A, $
we set \ $ (f + g)(x) = f(x) + g(x), \ (a f)(x) = a f(x), $
and $ (f g)(x) = f(x) \cdot g(x) \ (\forall x \in X). $ \
Then obviously $ \mathfrak{M}(X \rightarrow A) $ is a commutative
$ A-$algebra.
\par     \vskip  0.2 cm

{\bf Definition 4.10.} \ let $ (X, \widehat{d}_{X} ) $ be a
non-Archimedean metric space, $ (A, | \ |) $ be a
non-Archimedean valued ring, we define
$$ \widehat{BL}(X \rightarrow A) = \{ f \in
\mathfrak{M}(X \rightarrow A) : \ \| f \| _{\widehat{BL}}
< \infty \} \ \text{with} \ \| f \| _{\widehat{BL}} =
\max \{ \| f \| _{\infty }, \ \text{dil}(f) \},
$$ where $ \| f \| _{\infty } = \sup _{x \in X} | f(x)| $
is the supremum norm and $ \text{dil}(f) $ is the dilatation
of $ f $ as above.
\par     \vskip  0.2 cm

{\bf Theorem 4.11.} \ $ \widehat{BL}(X \rightarrow A) $ is an
$ A-$subalgebra of $ \mathfrak{M}(X \rightarrow A), $ and as an
$ A-$module, $ \| \cdot \| _{\widehat{BL}} $ is a
non-Archimedean norm on it.
Moreover, for any $ f, g \in \widehat{BL}(X \rightarrow A), $
we have $ \| f g \| _{\widehat{BL}} \leq \| f \| _{\widehat{BL}}
\cdot \| g \| _{\widehat{BL}}. $
\par     \vskip  0.1 cm
{\bf Proof.} \ For $ f, g \in \widehat{BL}(X \rightarrow A) $ and
$ a \in A, $ it follows easily by definition that
$  \| f + g \| _{\infty } \leq \max \{\| f \| _{\infty },
\| g \| _{\infty } \} < \infty $ and $  \| a f \| _{\infty } =
| a | \cdot  \| f \| _{\infty } < \infty . $ Also
\begin{align*}
&\text{dil}(f + g) = \sup _{x_{1}, x_{2} \in X \
\text{and} \ x_{1} \neq x_{2}} \frac{| (f(x_{1}) - f(x_{2})) +
(g(x_{1}) - g(x_{2}))|}{\widehat{d}_{X}(x_{1}, x_{2})} \\
&\leq \sup _{x_{1}, x_{2} \in X \
\text{and} \ x_{1} \neq x_{2}} \max \{\frac{| f(x_{1}) - f(x_{2})|}
{\widehat{d}_{X}(x_{1}, x_{2})}, \ \frac{| g(x_{1}) - g(x_{2})|}
{\widehat{d}_{X}(x_{1}, x_{2})} \} \\
&\leq \max \{\text{dil}(f), \text{dil}(g) \} < \infty ,
\quad \text{and} \\
&\text{dil}(a f) = \sup _{x_{1}, x_{2} \in X \
\text{and} \ x_{1} \neq x_{2}} \frac{| a f(x_{1}) - a f(x_{2})|}
{\widehat{d}_{X}(x_{1}, x_{2})} = | a | \cdot \text{dil}(f) < \infty .
\end{align*}
So $ \| f + g \| _{\widehat{BL}} < \infty $ and
$ \| a f \| _{\widehat{BL}} < \infty , $ and so $ f + g, \ a f \in
\widehat{BL}(X \rightarrow A). $ \\
From the above discussion, we also get $
\| f + g \| _{\widehat{BL}} \leq \max \{ \| f \| _{\widehat{BL}},
\| g \| _{\widehat{BL}} \} $ \ and \ $ \| a f \| _{\widehat{BL}} =
| a | \cdot \| f \| _{\widehat{BL}}. $
Moreover, it is obvious that $ \| f \| _{\widehat{BL}} \geq 0 $ and
$ \| f \| _{\widehat{BL}} = 0 $ if and only if $ f = 0. $ Therefore,
as an $ A-$module, $ \| \cdot \| _{\widehat{BL}} $ is a
non-Archimedean norm on $ \widehat{BL}(X \rightarrow A). $ \\
Now we come to show that $ \| f g \| _{\widehat{BL}} \leq
\| f \| _{\widehat{BL}} \cdot \| g \| _{\widehat{BL}}. $
In fact, since
\begin{align*}
&\| f g \| _{\infty } = \sup _{x \in X} | f(x)| \cdot | g(x)| \leq
\sup _{x \in X} | f(x)| \cdot \sup _{x \in X} | g(x)| =
\| f \| _{\infty } \cdot \| g \| _{\infty } \ \text{and} \\
&\text{dil}(f g) = \sup _{x_{1}, x_{2} \in X \
\text{and} \ x_{1} \neq x_{2}} \frac{| (f(x_{1}) - f(x_{2}))
g(x_{2}) + f(x_{1}) (g(x_{1}) - g(x_{2}))|}
{\widehat{d}_{X}(x_{1}, x_{2})} \\
&\leq \sup _{x_{1}, x_{2} \in X \
\text{and} \ x_{1} \neq x_{2}} \max \{\frac{| f(x_{1}) - f(x_{2})|
\cdot |g(x_{2})|}{\widehat{d}_{X}(x_{1}, x_{2})}, \
\frac{|f(x_{1})| \cdot | g(x_{1}) - g(x_{2})|}
{\widehat{d}_{X}(x_{1}, x_{2})} \} \\
&\leq \max \{\| g \| _{\infty } \cdot \text{dil}(f), \
\| f \| _{\infty } \cdot \text{dil}(g) \},
\end{align*}
it follows that $ \| f g \| _{\widehat{BL}} \leq
\max \{\| f \| _{\infty }, \text{dil}(f) \} \cdot
\max \{\| g \| _{\infty }, \text{dil}(g) \} =
\| f \| _{\widehat{BL}} \cdot \| g \| _{\widehat{BL}}. $
In particular, $ \| f g \| _{\widehat{BL}} < \infty , $
so $ f g \in \widehat{BL}(X \rightarrow A), $ which implies
that $ \widehat{BL}(X \rightarrow A) $ is
$ A-$subalgebra of $ \mathfrak{M}(X \rightarrow A), $ and the
proof is completed. \quad $ \Box $
\par     \vskip  0.2 cm

{\bf Remark 4.12.} \ For any non-Archimedean metric spaces
$ (X, \widehat{d}_{X} ) $ and $ (Y, \widehat{d}_{Y}), $ denote \\
$ \overline{\mathfrak{D}}_{\flat }(X, Y) =
\{\text{all maps} \ P : \overline{\mathcal{M}}_{\flat } (X)
\rightarrow Y \} = \mathfrak{M}(\overline{\mathcal{M}}_{\flat } (X)
\rightarrow Y). $ \\
Then one can similarly define and study non-Archimedean metrics
$ \widehat{\beta }_{X, Y}^{\lambda } $ and
$ \widehat{\beta }_{X, Y}^{\ast \lambda } $ for
$ \overline{\mathfrak{D}}_{\flat }(X, Y). $

\par     \vskip  0.3 cm

\hspace{-0.6cm}{\bf 5. \ Ultrametric structures on
non-Archimedean measures }

\par  \vskip  0.2 cm

Let $ K $ be a complete non-Archimedean valued field with absolute
value $ | \ |, $ and let $ (X, \widehat{d}_{X} ) $ be a compact
non-Archimedean metric space. It is well known that every non-empty
open subset in a non-Archimedean metric space is a disjoint union of
balls of the forms $ \overline{B}_{a} (r) $ (see [Sc], p.48). Then
it follows that an open subset of $ X $ is compact if and only if it
is a finite disjoint union of balls. In this section, we assume that
every ball in $ X $ and $ K $ contains at least two distinct points.
Let $ C(X) $ be the space of all continuous function $ f: X
\rightarrow K $ with the supremum norm $ \| f \|_{\infty } $ \ (As
sets, $ C(X) = \mathfrak{C}(X \rightarrow K)$). Recall that a $
K-$valued measure on $ X $ is a $ K-$linear functional $ \mu : \
C(X) \rightarrow K $ for which there is $ M \geq 0 $ such that the
inequality
$$ | \mu (f) | \leq M \| f \|_{\infty } $$ holds for every $ f \in C(X). $
Let $ \Omega (X) $ be the set of all open compact subsets of $ X. $
Let us denote by $ \chi _{A} $ the characteristic function of a set
$ A \subset X. $ Every measure $ \mu $ generates a mapping $ \mu :
\Omega (X) \rightarrow K $ by the rule $$ \mu (A): = \mu (\chi
_{A}), \quad A \in \Omega (X). $$ By this definition, $ \mu $ is
completely determined by its values on the balls, i.e., on the
elements of $ \mathcal{M}_{\flat }(X) $ (see, for example, Theorem
3.1 in $\S$5.3 of [Kh]). In particular, $ \mu $ can be considered as
a $ K-$valued function on $ \mathcal{M}_{\flat }(X), $ that is, the
restriction $ \mu _{| \mathcal{M}_{\flat }(X)} \in
\mathfrak{D}_{\flat }(X, K). $ \\
If $ \mu $ is a $ K-$valued measure on $ X, $ then for every $ f \in
C(X) $ we have
$$ \mu (f) = \int _{X} f(x) d \mu (x) = \lim _{\{x_{j} \in B_{j} \}}
\sum f(x_{j}) \mu (B_{j}), $$ where $ \sum f(x_{j}) \mu (B_{j}) $ is
the Riemann sums taken over increasingly fine covers of $ X $ by
mutually disjoint balls (see [MO], [MS], [Sch] and [Kh] for such
Monna-Springer integration theory). We denote
$$ \text{MEA}(X \rightarrow K) = \{\text{all} \
K-\text{valued measures on} \ X \}. $$ By Theorem 4.2.(1) above, it
is easy to see that all $ \widehat{\beta }_{X, K}^{\lambda } $ and $
\widehat{\beta }_{X, K}^{\ast \lambda } $ are non-Archimedean
metrics on $ \text{MEA}(X \rightarrow K), $ in particular, $
(\text{MEA}(X \rightarrow K), \widehat{\beta }_{X, K}) $
is a non-Archimedean metric space. Here \\
$ \widehat{\beta }_{X, K}^{\lambda }
(\mu ^{\prime }, \mu ^{\prime \prime}) = \widehat{\beta }_{X, K}^{\lambda }
(\mu ^{\prime } _{| \mathcal{M}_{\flat }(X)}, \
\mu ^{\prime \prime}_{| \mathcal{M}_{\flat }(X)}) $ and
$ \widehat{\beta }_{X, K}^{\ast \lambda }
(\mu ^{\prime }, \mu ^{\prime \prime}) =
\widehat{\beta }_{X, K}^{\ast \lambda }
(\mu ^{\prime } _{| \mathcal{M}_{\flat }(X)}, \
\mu ^{\prime \prime}_{| \mathcal{M}_{\flat }(X)}) \\
(\forall \mu ^{\prime }, \mu ^{\prime \prime} \in
\text{MEA}(X \rightarrow K)). $ \\
For a sequence of $ K-$valued measures $ \{\mu _{n} \}_{n =
1}^{\infty } $ and $ \mu $ in $ \text{MEA}(X \rightarrow K), $ it is
called that $ \mu _{n} $ converges to $ \mu , $ written $ \mu _{n}
\rightarrow \mu , $ if for every $ f \in \mathfrak{C}(X \rightarrow
K), \ \int _{X} f d \mu _{n} \rightarrow \int _{X} f d \mu $ as $ n
\rightarrow \infty , $ i.e., $ | \int _{X} f d \mu _{n} - \int _{X}
f d \mu | \rightarrow 0 $ as $ n \rightarrow \infty . $
\par     \vskip  0.2 cm

{\bf Definition 5.1.} \ Let $ (X, \widehat{d}_{X} ) $ be a
compact non-Archimedean metric space, and let $ K $ be a complete
non-Archimedean valued field with absolute value $ | \ |. $
For $ \mu _{1}, \mu _{2} \in \text{MEA}(X \rightarrow K) $ and
$ f \in \mathfrak{C}(X \rightarrow K), $ we denote
$ \int _{X} f d (\mu _{1} - \mu _{2}) = \int _{X} f d \mu _{1} -
\int _{X} f d \mu _{2}, $ and define the $ K-$valued Dudley metric (we will
show that it is indeed a metric) as follows:
$$ \widehat{D}_{X, K}(\mu _{1}, \mu _{2}) = \sup \{
| \int _{X} f d (\mu _{1} - \mu _{2})| : \ f \in
\mathfrak{C}(X \rightarrow K) \ \text{and} \ \| f \| _{\widehat{BL}}
\leq 1 \}. $$
(See Definition 4.10 above for $ \| f \| _{\widehat{BL}} $).
\par     \vskip  0.2 cm

{\bf Theorem 5.2.} \ Let $ (X, \widehat{d}_{X} ) $ be a
compact non-Archimedean metric space, and let $ K $ be a complete
non-Archimedean valued field with absolute value $ | \ |. $
Then $ \widehat{D}_{X, K} $ is a non-Archimedean metric
on $ \text{MEA}(X \rightarrow K). $
\par     \vskip  0.1 cm
{\bf Proof.} \ Obviously, $ \widehat{D}_{X, K}(\mu , \mu )
= 0 $ \ and \ $ \widehat{D}_{X, K}(\mu _{1}, \mu _{2}) =
\widehat{D}_{X, K}(\mu _{2}, \mu _{1}) \geq 0 $ for all $ \mu ,
\mu _{1}, \mu _{2} \in \text{MEA}(X \rightarrow K). $ Now let
$ \mu _{1}, \mu _{2}, \mu _{3} \in \text{MEA}(X \rightarrow K), $
since
\begin{align*}
&| \int _{X} f d (\mu _{1} - \mu _{2})| =
| \int _{X} f d (\mu _{1} - \mu _{3}) +
\int _{X} f d (\mu _{3} - \mu _{2})| \\
&\leq \max \{
| \int _{X} f d (\mu _{1} - \mu _{3})|, \
| \int _{X} f d (\mu _{3} - \mu _{2})| \} \
(\forall f \in \mathfrak{C}(X \rightarrow K) \ \text{with} \
\| f \| _{\widehat{BL}} \leq 1), \\
&\text{we have} \ \widehat{D}_{X, K}(\mu _{1}, \mu _{2}) \leq \\
&\sup \{ \max \{ | \int _{X} f d (\mu _{1} - \mu _{3})|,
| \int _{X} f d (\mu _{3} - \mu _{2})| \} : \
f \in \mathfrak{C}(X \rightarrow K) \ \text{with} \
\| f \| _{\widehat{BL}} \leq 1 \} \\
&= \max \{\widehat{D}_{X, K}(\mu _{1}, \mu _{3}),
\widehat{D}_{X, K}(\mu _{3}, \mu _{2}) \}.
\end{align*}
Hence the strong triangle inequality holds. \\
Lastly, we assume that $ \widehat{D}_{X, K}(\mu _{1}, \mu _{2})
= 0, $ and we need to show that $ \mu _{1} = \mu _{2}. $ To see this,
by the assumption, we have $ | \int _{X} f d (\mu _{1} - \mu _{2})|
= 0, $ and hence \\
$ \int _{X} f d \mu _{1} =
\int _{X} f d \mu _{2} $ for all $ f \in \mathfrak{C}(X \rightarrow K) $
with $ \| f \| _{\widehat{BL}} \leq 1. \quad (\star ) $ \\
Take $ B \in \mathcal{M}_{\flat }(X), $ and let
$ \chi _{B} $ be the characteristic function, i.e., \\
$ \chi _{B}(x) = \left\{ \begin{array}{l} 1 \quad \text{if}
\ x \in B, \\
0 \quad \text{otherwise.}
\end{array} \right. $
Then $ \chi _{B} $ is a locally constant function, hence a continuous
function on $ X $ (see [K], p.31), i.e.,
$ \chi _{B} \in \mathfrak{C}(X \rightarrow K). $ Moreover,
it follows easily by
definition that $ \| \chi _{B} \| _{\infty } = 1 $ and
$ \text{dil}(\chi _{B}) \leq 1 / \text{diam} B. $ If $ \text{diam} B
\geq 1, $ then $ \text{dil}(\chi _{B}) \leq 1, $ so
$ \| \chi _{B} \| _{\widehat{BL}} = \| \chi _{B} \| _{\infty } = 1, $
hence by the above $ (\star ), $ we get $ \mu _{1}(B) =
\int _{X} \chi _{B} d \mu _{1} = \int _{X} \chi _{B} d \mu _{2} =
\mu _{2}(B). $ If $ \text{diam} B < 1, $ then note that $ \text{diam} B > 0, $
we can take an element $ c \in K $ such that $ 0 < | c | \leq
\text{diam} B < 1. $ Let $ f = c \cdot \chi _{B}, $ then $ f \in
\mathfrak{C}(X \rightarrow K), \ \| f \| _{\infty } =
| c | \cdot \| \chi _{B} \| _{\infty } = | c |, $ and
$ \text{dil}(f) = | c | \cdot \text{dil}(\chi _{B}) \leq
| c | / \text{diam} B \leq 1, $ which implies
$ \| f \| _{\widehat{BL}} \leq 1. $ Then by the above $ (\star ), $
we have $ \int _{X} f d \mu _{1} = \int _{X} f d \mu _{2}, $
i.e., $ \int _{X} c \cdot \chi _{B} d \mu _{1} =
\int _{X} c \cdot \chi _{B} d \mu _{2}, $ so
$ \int _{X} \chi _{B} d \mu _{1} = \int _{X} \chi _{B} d \mu _{2}, $
i.e., $ \mu _{1}(B) = \mu _{2}(B). $ To sum up, we have shown that
$ \mu _{1}(B) = \mu _{2}(B) $ for every $ B \in \mathcal{M}_{\flat }(X), $
which implies $ \mu _{1} = \mu _{2}. $ Therefore,
$ \widehat{D}_{X, K} $ is a non-Archimedean metric
on $ \text{MEA}(X \rightarrow K). $ The proof is completed.
\quad $ \Box $
\par     \vskip  0.2 cm
For $ X $ and $ K $ in Definition 5.1, recall that \\
$ \widehat{BL}(X \rightarrow K) = \{ f \in
\mathfrak{M}(X \rightarrow K) : \ \| f \| _{\widehat{BL}}
< \infty \}. $
\par     \vskip  0.2 cm

{\bf Theorem 5.3.} \ Let $ (X, \widehat{d}_{X} ) $ be a
compact non-Archimedean metric space, and let $ K $ be a complete
non-Archimedean valued field with absolute value $ | \ |. $
For $ \mu $ and a sequence $ \{\mu _{n} \}_{n = 1}^{\infty } $ in
$ \text{MEA}(X \rightarrow K), $ if there exists a $ c > 0 $ such that
$ \| \mu _{n} \| \leq c $ for all positive integers $ n, $
then the following statements are equivalent: \\
(1) \ $ \mu _{n} \rightarrow \mu . $ \\
(2) \ $ \int _{X} f d \mu _{n} \rightarrow \int _{X} f d \mu $ for
all $ f \in \widehat{BL}(X \rightarrow K). $ \\
(3) \ $ \widehat{D}_{X, K}(\mu _{n}, \mu ) \rightarrow 0. $
\par     \vskip  0.1 cm
{\bf Proof.} \ (1) $ \Longrightarrow $ (2). Let
$ f \in \widehat{BL}(X \rightarrow K), $ then $ \| f \| _{\widehat{BL}}
< \infty , $ in particular, $ \text{dil}(f) < \infty . $ So $ f $ is
uniformly continuous, and $ f \in \mathfrak{C}_{u}(X \rightarrow K)
\subset \mathfrak{C}(X \rightarrow K), $ hence
$ \int _{X} f d \mu _{n} \rightarrow \int _{X} f d \mu . $ \\
(2) $ \Longrightarrow $ (1). Let $ \text{Step}(X) = \{K-\text{valued
locally constant functions on} \ X \}, $ then $ \text{Step}(X) $ is
dense in the $ K-$Banach space $ \mathfrak{C}(X \rightarrow K) $ with
the supremum norm $ \| \cdot \| _{\infty } $ defined above
(see [Wa], p.237). We assert that $ \text{Step}(X)
\subset \widehat{BL}(X \rightarrow K). $ In fact, since $ X $ is
compact, for every $ g \in \text{Step}(X), g $ is
a finite linear combination of characteristic functions of disjoint
balls, i.e., $ g = \sum _{1 \leq i \leq r} a_{i} \chi _{i}, $ where
$ a_{i} \in K, \chi _{i} $ is the characteristic function of the
ball $ B_{i} $ and $ X = \sqcup _{1 \leq i \leq r} B_{i} \ $
(the disjoint union). Obviously, $ \| g \| _{\infty } =
\max \{ | a_{i} | : \ i = 1, \cdots , r \} < \infty . $
Fix a point $ x_{i} \in B_{i} $ for each $ i, $
then $ \text{dist}(B_{i}, B_{j}) = \widehat{d}_{X}(x_{i}, x_{j}) $
for each pair $ (i, j) $ (see [Sc], p.48). Denote
$ \gamma = \min \{\widehat{d}_{X}(x_{i}, x_{j}) : \
1 \leq i, j \leq r \ \text{and} \ i \neq j \} $
and $ \alpha = \max \{|a_{i} -a_{j}| : \
1 \leq i, j \leq r \}, $ then $ \gamma > 0, 0 \leq \alpha < \infty , $
and it is easy to verify that
$ \text{dil}(g) \leq \alpha / \gamma < \infty , $
hence $ \| g \| _{\widehat{BL}} < \infty . $ Therefore
$ \text{Step}(X) \subset \widehat{BL}(X \rightarrow K). $ \\
Now for every $ f \in \mathfrak{C}(X \rightarrow K), $ since
$ \text{Step}(X) $ is dense in $ \mathfrak{C}(X \rightarrow K), $
there exists a Cauchy sequence $ \{f _{m} \}_{m = 1}^{\infty } $ in
$ \text{Step}(X) $ such that $ f _{m} \rightarrow f. $ Then by the
proof of Proposition 12.1 in [Wa], we have
$$ \lim _{m \rightarrow \infty } \int _{X} f _{m} d \mu =
\int _{X} f d \mu \ \text{and} \
\lim _{m \rightarrow \infty } \int _{X} f _{m} d \mu _{n} =
\int _{X} f d \mu _{n} \ (n =1, 2, \cdots ). $$
Moreover, for every such $ f _{m}, $ as shown in the above assertion,
it is also in $ \widehat{BL}(X \rightarrow K), $ so by the assumption,
we have $ \int _{X} f _{m} d \mu _{n} \rightarrow
\int _{X} f _{m} d \mu $ as $ n \rightarrow \infty . $ \\
Denote $ c_{0} = \max \{\| \mu \|, \ c \}. $ Then by
the hypothesis, $ \| \mu _{n} \| \leq c_{0} $ for all
positive integers $ n. $ Now for any $ \varepsilon > 0, $ since
$ f _{m} \rightarrow f $ under the supremum norm $ \| \cdot \| _{\infty }, $
there exists a $ M > 0 $ such that for all integers $ m > M $ we have
$ \| f _{m} - f \| _{\infty } < \varepsilon / c_{0}, $ in particular,
$ | f _{m}(x) - f(x)| < \varepsilon / c_{0} $ for all $ x \in X. $
So $ | \int _{X} f _{m} d \mu _{n} - \int _{X} f d \mu _{n}| <
 \varepsilon / c_{0} \cdot c_{0} = \varepsilon $ for all positive
integers $ n $ and all integers $ m > M $ (see [K], p.40). Also
$ | \int _{X} f _{m} d \mu - \int _{X} f d \mu | < \varepsilon $
for all integers $  m > M. $ We fix an integer $ m_{0} > M, $ then from
$ \int _{X} f _{m_{0}} d \mu _{n} \rightarrow
\int _{X} f _{m_{0}} d \mu $ as $ n \rightarrow \infty , $ there exists
a $ N > 0 $ such that for all integers $ n > N, $ we have
$ | \int _{X} f _{m_{0}} d \mu _{n} - \int _{X} f _{m_{0}} d \mu |
< \varepsilon , $ and then
\begin{align*}
&| \int _{X} f d \mu _{n} - \int _{X} f d \mu | \leq \\
&\max \{| \int _{X} f d \mu _{n} - \int _{X} f _{m_{0}} d \mu _{n}|, \
| \int _{X} f _{m_{0}} d \mu _{n} - \int _{X} f _{m_{0}} d \mu |, \
| \int _{X} f _{m_{0}} d \mu - \int _{X} f d \mu | \} \\
&< \varepsilon .
\end{align*}
So $ \int _{X} f d \mu _{n} \rightarrow \int _{X} f d \mu $ as $ n
\rightarrow \infty , $ hence $ \mu _{n} \rightarrow \mu . $ \\
(3) \ $ \Longrightarrow $ (1). \ Let $ B \in \mathcal{M}_{\flat }(X)
$ and $ \chi _{B} $ be the corresponding characteristic function.
Take $ a \in K $ such that $ 0 < | a | \leq \text{diam}(B). $
We set $ \chi = \left\{\begin{array}{l} \chi _{B} \quad \text{if}
\ \text{diam}(B) \geq 1, \\
a \cdot \chi _{B} \quad \text{if} \ \text{diam}(B) < 1.
\end{array} \right. $
Then from the proof of the above Theorem 5.2, we know that $ \| \chi
\| _{\widehat{BL}} \leq 1, $ and then $ | \int _{X} \chi d \mu _{n}
- \int _{X} \chi d \mu | \leq \widehat{D}_{X, K}(\mu _{n}, \mu )
\rightarrow 0, $ so $ \int _{X} \chi d \mu _{n} \rightarrow \int
_{X} \chi d \mu , $ and so $ \int _{X} \chi _{B} d \mu _{n}
\rightarrow \int _{X} \chi _{B} d \mu . $ It then follows that $
\int _{X} g d \mu _{n} \rightarrow \int _{X} g d \mu $ for
all $ g \in \text{Step}(X). $ \\
Now let $ f \in \mathfrak{C}(X \rightarrow K), $ similar to the
proof of the above $ ' (2)\Rightarrow (1)', $ there exists a Cauchy
sequence $ \{f_{m} \}_{m = 1}^{\infty } $ in $ \text{Step}(X), $
such that $ f_{m} \rightarrow f. $ By the above discussion, we have
$ \int _{X} f_{m} d \mu _{n} \rightarrow \int _{X} f_{m} d \mu \ (m
=1, 2, \cdots ). $ Denote $ c_{0} = \max \{ \| \mu \|, c \}. $ For
any $ \varepsilon > 0, $ since $ f_{m} \rightarrow f, $ as discussed
in the proof of the above $ '(2)\Rightarrow (1)', $ there exists a
positive number $ M $ such that for all integers $ m
> M, $ we have $ | \int _{X} f_{m} d \mu _{n} - \int _{X} f d \mu
_{n}| < \varepsilon $ for all integers $ n > 0, $ and $ | \int _{X}
f_{m} d \mu - \int _{X} f d \mu| < \varepsilon . $ Fix an integer $
m_{0} > M, $ then since $ \int _{X} f_{m_{0}} d \mu _{n} \rightarrow
\int _{X} f_{m_{0}} d \mu $ as $ n \rightarrow \infty , $ there
exists a $ N > 0 $ such that for all $ n > N, $ we have $ | \int
_{X} f_{m_{0}} d \mu _{n} - \int _{X} f_{m_{0}} d \mu | <
\varepsilon , $ then as done in the proof of the above $ '
(2)\Rightarrow (1)', $ we have $ | \int _{X} f d \mu _{n} - \int
_{X} f d \mu | < \varepsilon . $ Hence $ \mu _{n} \rightarrow \mu .
$ \\
(1) \ $ \Longrightarrow $ (3). \ Write $ \mathfrak{F} = \{ f \in
\mathfrak{C}(X \rightarrow K): \ \| f \| _{\widehat{BL}} \leq 1 \}.
$ By definition, for any $ f \in \mathfrak{F}, $ we have $
\text{dil}(f) \leq \| f \| _{\widehat{BL}} \leq 1, $ so $ \mid
f(x_{1}) - f(x_{2}) \mid \leq \text{dil}(f) \cdot
\widehat{d}_{X}(x_{i}, x_{j}) \leq \widehat{d}_{X}(x_{i}, x_{j}) \
(\forall x_{1}, x_{2} \in X ). $
Hence for any $ \varepsilon > 0, $ we have that
$ \widehat{d}_{X}(x_{i}, x_{j}) < \varepsilon $ implies
$ \mid f(x_{1}) - f(x_{2}) \mid < \varepsilon $ for all
$ x_{1}, x_{2} \in X $ and all $ f \in
\mathfrak{F}, $ i.e., $ \mathfrak{F} $ is uniformly equicontinuous
(see [D], p.51), in particular, $ \mathfrak{F} $ is equicontinuous.
Also we have $ \| f \| _{\infty } \leq \| f \| _{\widehat{BL}} \leq
1 $ for all $ f \in \mathfrak{F}. $ So $ \mathfrak{F} $ is uniformly
bounded, hence by the ultrametric version of the Ascoli Theorem (see
Thm.3.8.2 in [PS], p.144), we know that $ \mathfrak{F} $ is
compactoid in $ \mathfrak{C}(X \rightarrow K). $ Therefore, for any
$ \varepsilon > 0, $ there exists a finite set $ f_{1}, \cdots ,
f_{r} \in \mathfrak{C}(X \rightarrow K) $ such that for every $ f
\in \mathfrak{F}, \ \| f - f_{j} \| _{\infty } < \varepsilon $ for
some $ j. $ Since $ \mu _{n} \rightarrow \mu , $ there exists a
positive integer $ N $ such that for every integer $ n > N, $ we
have $ \mid \int _{X} f_{i} d\mu _{n} - \int _{X} f_{i} d\mu \mid <
\varepsilon $ for each $ i = 1, \cdots , r. $ Thus,
\begin{align*}
&\mid \int _{X} f d( \mu _{n} - \mu )\mid = \mid \int _{X} f d \mu
_{n} - \int _{X} f d \mu \mid \\
&= \mid \int _{X} (f - f_{j}) d \mu _{n} - \int _{X} (f - f_{j}) d
\mu + \int _{X} f_{j} d \mu _{n} - \int _{X} f_{j} d \mu \mid \\
&\leq \max \{\mid \int _{X} (f - f_{j}) d \mu _{n} \mid, \mid \int
_{X} (f - f_{j}) d \mu \mid, \mid \int _{X} f_{j} d \mu _{n} - \int
_{X} f_{j} d \mu \mid \}.
\end{align*}
Since $ \| f - f_{j} \| _{\infty } < \varepsilon $ and $ \| \mu _{n}
\| \leq c_{0}, \| \mu \| \leq c_{0} \ (c_{0}= \max \{c, \| \mu \|
\}), $ we have $ \mid \int _{X} (f - f_{j}) d \mu _{n} \mid \leq
c_{0} \cdot \varepsilon  $ and $ \mid \int _{X} (f - f_{j}) d \mu
\mid \leq c_{0} \cdot \varepsilon $ (see [K], p.40). Therefore, $
\mid \int _{X} f d( \mu _{n} - \mu )\mid \leq \max \{c_{0} \cdot
\varepsilon , \varepsilon \} = \max \{1, c_{0} \} \cdot \varepsilon
\ (\forall f \in \mathfrak{F}). $ Hence by definition, $
\widehat{D}_{X, K}(\mu _{n}, \mu ) \leq \max \{1, c_{0} \} \cdot
\varepsilon , $ which implies $ \widehat{D}_{X, K}(\mu _{n}, \mu )
\rightarrow 0. $ The proof is completed. \quad $ \Box $
\par     \vskip  0.2 cm

{\bf Remark 5.4.} \ Let $ X, K, \mu $ and $ \{\mu _{n} \}_{n =
1}^{\infty } $ be as in Theorem 5.3 above,
then we also have \\
 $ \widehat{\beta }_{X, K}(\mu _{n}, \mu ) \rightarrow 0
\Longrightarrow \mu _{n} \rightarrow \mu . $
\par     \vskip  0.1 cm
In fact, if $ \widehat{\beta }_{X, K}(\mu _{n}, \mu ) \rightarrow 0,
$ then by Theorem 4.2.(2) above, we have $ \widehat{\rho }_{s}(\mu
_{n}, \mu ) \rightarrow 0, $ i.e., $ \sup _{B \in \mathcal{M}_{\flat
}(X)} | \mu _{n}(B) - \ \mu (B)| \rightarrow 0, \ \text{and then} \
| \int _{X} \chi _{B} d \mu _{n} - \int _{X} \chi _{B} d \mu | = |
\mu _{n}(B) - \mu(B)| \rightarrow 0, $ i.e., $ \int _{X} \chi _{B} d
\mu _{n} \rightarrow \int _{X} \chi _{B} d \mu , $ where $ \chi _{B}
$ is the characteristic function of $ B \ ( \forall B \in
\mathcal{M}_{\flat }(X)). $ It follows then $ \int _{X} g d \mu _{n}
\rightarrow \int _{X} g d \mu $ for all $ g \in \text{Step}(X), $
and as done in the proof of the above Theorem 5.3, one can similarly
show that, for every $ f \in \mathfrak{C}(X \rightarrow K), \ \int
_{X} f d \mu _{n} \rightarrow \int _{X} f d \mu $ as $ n \rightarrow
\infty , $ so $ \mu _{n} \rightarrow \mu . $  \quad $ \Box $
\par     \vskip  0.1 cm
A question here is that is it also true that $ \mu _{n} \rightarrow
\mu \Longrightarrow \widehat{\beta }_{X, K}(\mu _{n}, \mu )
\rightarrow 0 ? $

\par  \vskip 0.3 cm

{ \bf Acknowledgments. } \ I would like to thank the anonymous
referee for a very careful reading of the paper and many helpful
comments, especially for pointing out the use of the Monna-Springer
integration in Section 5 of the paper.

\par  \vskip 0.2 cm

\hspace{-0.8cm} {\bf References }
\begin{description}

\item[[BBI]] D. Burago, Y. Burago, S. Ivanov, A Course in Metric
Geometry, Providence, Rhode Island: American Mathematical Society,
2001.

\item[[BGR]] S. Bosch, U. Guntzer, R. Remmert, Non-Archimedean
Analysis, Berlin: Springer-Verlag, 1984.

\item[[D]] R. M. Dudley, Real Analysis and Probability,
Second Edition, Cambridge: Cambridge University Press, 2003.

\item[[FP]] J. Frensel, M. Van der Put, Rigid Analytic Geometry and
Its Applications, Progress in Math., Vol. 218, Boston: Birkhauser, 2004.
\item[[G]] M. Gromov, Metric Structures for Riemannian and
Non-Riemannian Spaces, Boston: Birkhauser, 2001.

\item[[K]] N. Koblitz, $ p-$adic Numbers, $ p-$adic Analysis, and
Zeta Functions, Second Edition, New York: Springer-Verlag, 1984.

\item[[Kh]] A. Yu. Khrennikov, Non-Archimedean analysis and its applications,
Moscow, Fiziko-Mathematichaskay Literatura, 2009 (in Russian).

\item[[L]] S. Lang, Algebraic Number Theory, Second Edition, New
York: Springer-Verlag, 1994.

\item[[M]] J.R. Munkres, Topology, 2nd Edition, Beijing: China
Machine Press, 2004.

\item[[Mo]] A. Monna, Analyse non-Archimedean, Springer, New York, 1970.

\item[[MS]] A. Monna, T. Springer, Integration non-Archimedean I-II,
Indag. Math., 25(1963), no.4, 634-653.

\item[[PS]] C. Perez-Garcia and W. H. Schikhof, Locally Convex
Spaces over Non-Archimedean Valued Fields, Cambridge: Cambridge
University Press, 2010.

\item[[Q]] D. R. Qiu, Geometry of non-Archimedean Gromov-Hausdorff
distance, $ p-$adic numbers, ultrametric analysis and applications,
4 (2009), 317-337.

\item[[Ra]] S. Rachev, Hausdorff metric structure of the space of
probability measures, Translated from Zapiski Nauchnykh Seminarov
Leningradskogo Otdeleniya Matematicheskogo Instituta im. V. A.
Steklova AN SSSR, 87 (1979), 87-103.

\item[[Sc]] W. H. Schikhof, Ultrametric Calculus, London:
Cambridge University Press, 1984.

\item[[Sch]] W. H. Schikhof, Non-Archimedean Harmonic Analysis, Catholic Univ. Press,
Nijmegen, 1967.

\item[[Se]] J.-P. Serre, Local Fields, New York: Springer-Verlag,
1979.

\item[[Wa]] L.C. Washington, Introduction to Cyclotomic Fields, 2nd Edition,
New York: Springer-Verlag, 1997.

\end{description}

\end{document}